\documentclass[11pt]{amsart}
\usepackage{amssymb}
\usepackage{amsfonts}
\usepackage{color}

\input{diagrams}

\newarrow{Dash}{}{dash}{}{dash}>
\newarrow{Into}>--->
\newarrow{Sub}C--->
\newarrow{Onto}----{>>}
\newarrow{Equal}===={=}

\newtheorem{theorem}{Theorem}[section]
\newtheorem{proposition}[theorem]{Proposition}
\newtheorem{lemma}[theorem]{Lemma}
\newtheorem{corollary}[theorem]{Corollary}
\newtheorem{definition}[theorem]{Definition}

\newtheorem{example}[theorem]{Example}

\theoremstyle{remark}
\newtheorem{remark}[theorem]{Remark}

\renewenvironment{proof}{{\noindent\bf Proof.}}{\hfill $\Box$\par\vskip3mm}

\newcommand{\Nat}{\mathbf\underline{\underline{Nat}}\,}
\newcommand{\Ker}{{\rm Ker}\,}

\newcommand{\im}{{\rm Im}\,}

\newcommand{\Hom}{{\rm Hom}}

\newcommand{\Cc}{\mathcal{C}}
\newcommand{\Dd}{\mathcal{D}}

\newcommand{\Mm}{\mathcal{M}}

\def\NN{{\mathbb N}}

\def\text#1{{\rm {\rm #1}}}

\def\dul#1{\underline{\underline{#1}}}
\def\Nat{\dul{\rm Nat}}
\def\Id{\mathbb{I}}
\def\*C{{{}^*\mathcal C}}

\begin{document}
\title{Frobenius Extensions of Corings}

\begin{abstract}
Let $\Cc$ and $\Dd$ be two corings over a ring $A$ and $\Cc\stackrel{\lambda}{\longrightarrow}\Dd$ be a morphism of corings. We investigate the situation when the associated induced ("corestriction of scalars") functor $\Mm^\Cc\longrightarrow \Mm^\Dd$ is a Frobenius functor, and call these morphisms Frobenius extensions of corings. The characterization theorem generalizes notions such as Frobenius corings and is applied to several situations; in particular, provided some (general enough) flatness conditions hold, the notion proves to be dual to that of Frobenius extensions of rings (algebras). Several finiteness theorems are given for each case we consider; these theorems extend existing results from Frobenius extensions of rings or from Frobenius corings, showing that a certain finiteness property almost always occur for many instances of Frobenius functors.
\end{abstract}

\author{Miodrag Cristian Iovanov}
\thanks{$^*$ This paper was partially supported by a CNCSIS BD-type grant and was written within the frame of the bilateral Flemish-Romanian project "New Techniques in Hopf Algebra Theory and Graded Rings"}
\date{}
\maketitle

\section*{Introduction}
Let $A$ be an algebra over a field $K$. $A$ is called Frobenius if the right regular $A$-module $A$ is isomorphic to the $K$-dual of $A$, $A^*$. More generally, if $\varphi:A\longrightarrow B$ is a morphism of rings, then this is called a Frobenius extension if $B$ is finitely generated projective as left $A$-module and $B\simeq \Hom_A(B,A)$. Equivalent characterizations of this concept have lead to the introduction of the concept of Frobenius functor: a functor $F:\Cc\longrightarrow \Dd$ between two categories $\Cc$ and $\Dd$ is called Frobenius if and only if it has both left and right adjoints and its left adjoint is naturally isomorphic to its right adjoint. This concept was first introduced by Morita in \cite{Mo}, and is inspired by an equivalent characterization of Frobenius extensions of rings: an extension $A\stackrel{\varphi}{\longrightarrow} B$ is Frobenius if and only if the associated forgetful functor from $\Mm_B$ to $\Mm_A$ is Frobenius in the above specified sense. It turns out that this concept is left right symmetric. Further notions have been introduced and studied in connection to the notion of Frobenius functor. \\
If $\Cc$ is a coring over a ring $A$, with counit $\varepsilon_\Cc:\Cc\longrightarrow A$, then there is an associated induced (forgetful) functor $U$ going from the category $\Mm^\Cc$ of right comodules over $\Cc$ to $\Mm_A$ associating to each comodule $M$ the underlying $A$-module $M$. $\Cc$ is called a Frobenius coring provided that $U$ is a Frobenius functor. Given an extension of rings $\varphi:A\longrightarrow B$, a certain $B$-coring structure can be introduced on $B\otimes_AB$ called canonical Sweedler coring structure which connects the notion of Frobenius extension of rings to that of Frobenius corings, i.e. given certain flatness conditions hold, $A\stackrel{\varphi}{\longrightarrow} B$ is a Frobenius extension if and only if the $B$-coring $B\otimes_AB$ is Frobenius. \\
In this paper we aim to introduce and study the concept of Frobenius extensions of Corings, a notion which generalizes all these notions and also recovers some other known instances of Frobenius functors, as for example, the finite coproduct functor for a category of modules (or more generally, of comodules over a coring). Given two $A$-corings $\Cc$ and $\Dd$ and a morphism of corings $\lambda:\Cc\longrightarrow \Dd$, we can associate a "forgetful functor" (called corestriction functor) $U$ from $\Mm^\Cc$ to $\Mm^\Dd$, which associates to each $\Cc$-comodule $(M,\rho_M)$ the $\Dd$-comodule $(M,(M\otimes \lambda)\circ\rho_M)$, where the identity $1_M$ of a module $M$ will be often denoted simply by $M$. The extension $\Cc\stackrel{\lambda}{\longrightarrow}\Dd$ will be called Frobenius if $U$ is a Frobenius functor. When certain flatness conditions hold, $U$ has a right adjoint $F$ constructed by using the cotensor product: $F=-\square_\Dd\Cc$. The flatness conditions required will prove not to be very restrictive allowing us to apply the considerations to several situations, and therefore look for equivalent conditions for $F$ to be also a left adjoint for $U$. The characterization formula we obtain will prove to be a close dualization of the one characterizing the case of Frobenius extensions of rings, and in particular for the case of extensions of coalgebras we obtain a theorem which is dual to the one characterizing Frobenius extensions of algebras. \\
It is notable that in the case of Frobenius extensions of rings as well as for Frobenius corings certain finiteness theorems hold: if the functor $U$ associated to an extension of rings $A\stackrel{\varphi}{\longrightarrow} B$ (or to an $A$-coring $\Cc$) is Frobenius then $B$ is finitely generated projective as left (and as right) $A$-module (or $\Cc$ is finitely generated projective as left and as right $A$-module). We aim to prove several such finiteness theorems. For each situation of Frobenius extension of corings we consider we give a finiteness theorem. One application is for the coproduct of comodules indexed by a set $I$ functor; it is known from \cite{I1} that the coproduct is a Frobenius functor (equivalently, the product of comodules exist and is isomorphic to the coproduct) if and only if the index set $I$ is finite; this finiteness property will follow also as an application of the characterization of Frobenius extensions of corings Theorem. Another application is with Frobenius extensions of coalgebras. We prove that if $C\stackrel{\lambda}{\longrightarrow}D$ is a Frobenius extension of coalgebras over a field then ${\rm dim}(C)\leq {\rm dim}(D)$ or they are both finite. Finally, we give some connections between Frobenius extensions of coalgebras $C\stackrel{\lambda}{\longrightarrow}D$ and the Frobenius property of the dual extension of algebras $D^*\stackrel{\lambda^*}{\longrightarrow}{C^*}$, showing that these are equivalent whenever $C$ or $D$ is finite dimensional. 

\section{Extensions of Corings and associated functors}
Let $A$ be a ring and $\Cc$ and $\Dd$ be two $A$-corings. Denote by $\Delta_\Cc$, $\varepsilon_\Cc$ and $\Delta_\Dd$, $\varepsilon_\Dd$ the comultiplication and respectively counit of $\Cc$ and $\Dd$ respectively. Our convention is to write $\otimes$ every time we write a tensor product over $A$ of $A$-modules. For $c\in \Cc$ we use the Sweedler notation $\Delta_\Cc(c)=c_1\otimes c_2\in \Cc\otimes \Cc$ with the omitted summation symbol. Also, if $M$ is a right (or left) $\Cc$ comodule with comultiplication $\rho:M\rightarrow M\otimes \Cc$ (or $\rho:M\rightarrow \Cc\otimes M$) we write $\rho(m)=m_0\otimes m_1$ (or $m_{-1}\otimes m_0$).  For basic facts on corings and their comodules the reader is referred to \cite{BW}. Let $(\Cc,\Delta_\Cc,\varepsilon_\Cc)$ and $(\Dd,\Delta_\Dd,\varepsilon_{\Dd})$ be $A$-corings and $\lambda:\Cc\rightarrow \Dd$ be a morphism of corings, that is, $\lambda(c)_1\otimes\lambda(c)_2=\lambda(c_1)\otimes\lambda(c_2)$ and $\varepsilon_D(\lambda(c))=\varepsilon_C(c)$. Then there is an associated functor $U:\Mm^{\Cc}\rightarrow \Mm^\Dd$, called the corestriction functor and defined by $U(M)=M$ with the right $\Dd$-comodule structure given by $m\mapsto m_0\otimes\lambda(m_1)$. In particular, $\Cc$ has a $\Dd$-bicomodule structure. Recall that for a right $\Cc$-comodule $M$ and a left $\Cc$-comodule $N$ the cotensor product of $M$ and $N$ over $\Cc$, $M\square_\Cc N$ is defined as follows: let $\omega_{M,N}:M\otimes N\rightarrow M\otimes \Cc\otimes N$, $\omega_{M,N}=\rho_M\otimes\Cc-\Cc\otimes\rho_N$. Then $M\square_\Cc N$ is defined as the abelian group arising as the kernel of $\omega_{M,N}$, that is we have an exact sequence
\begin{diagram}
0 & \rTo & M\square_\Cc N & \rTo M & \otimes N & \rTo^{\omega_{M,N}} & M\otimes \Cc\otimes N
\end{diagram}
If $f:M\longrightarrow M'$ is a morphism of right $\Cc$-comodules and $g:N\longrightarrow N'$ is a morphism of left $\Cc$-comodules then one can define $f\square_{\Cc}g:M\square_{\Cc}N\longrightarrow M'\square_{\Cc}N'$ by $(f\square_{\Cc}g)(m\otimes n)=f(m)\otimes g(n)$. 
Let $N$ be a right $\Dd$-comodule. If ${}_A\Cc$ is flat, then by tensoring the exact sequence
$$0\longrightarrow N\square_\Dd \Cc \longrightarrow N\otimes\Cc \stackrel{\omega_{N,\Cc}}{\longrightarrow} N\otimes \Dd\otimes \Cc$$
with $\Cc$ on the right hand side we obtain the exact sequence
$$0\longrightarrow (N\square_\Dd \Cc)\otimes\Cc \longrightarrow N\otimes\Cc\otimes \Cc \stackrel{\omega_{N,\Cc}\otimes\Cc}{\longrightarrow} N\otimes \Dd\otimes \Cc$$
and therefore $(N\square_\Dd\Cc)\otimes\Cc\cong N\square_\Dd(\Cc\otimes\Cc)$ by the natural mapping $n\otimes c\otimes c'\mapsto n\otimes c\otimes c'$. This allows us to define a right $\Cc$-comodule structure on $N\square_\Dd\Cc$ by $N\square_\Dd \Cc\ni n\otimes c\mapsto n\otimes c_1\otimes c_2\in N\square_\Dd(\Cc\otimes\Cc)\cong(N\square_\Dd\Cc)\otimes\Cc$. More generally, we say that the coring morphism $\lambda$ is (right) {\it pure} if for every right $\Cc$-comodule $N$ the morphism $\omega_{N,C}$ is $\Cc$-{\it pure} in $\Mm_{A}$ (see \cite{BW}, Section 24, 24.8). In this case, the above isomorphism always holds.

Then we can define a functor $F$ from $\Mm^\Dd$, the category of all right $\Dd$-comodules to $\Mm^\Cc$ by writing $F(N)=N\square_\Dd \Cc$. If $f:N\rightarrow N'$ is a morphism in $\Mm^\Dd$, then $F(f)=f\square_\Dd C$, $F(f)(n\otimes c)=f(n)\otimes c$. By \cite{BW}, Section 22 (22.10) and Section 24 (24.11) we have that if ${}_A\Cc$ is flat or more generally, if $\lambda$ is a pure morphism of corings, then $F$ is right adjoint to $U$. Therefore, in this case the functor $U$ is Frobenius if and only if $F$ is also a left adjoint to $U$, because any two left (or right) adjoints of a functor are naturally equivalent (Kan, \cite{Kn}). In order even have a functor $F:\Mm^\Dd\longrightarrow \Mm^\Cc$, we will always assume that either $\lambda$ is a pure morphism of corings or ${}_A\Cc$ is flat or, more generally, the canonical morphism $(N\square_\Dd\Cc)\otimes\Cc\stackrel{\sim}{\rightarrow} N\square_\Dd(\Cc\otimes\Cc)$ is an isomorphism (equivalently, $\omega_{N,\Cc}$ is $\Cc$-pure in $\Mm_A$; see \cite{BW} 21.4). Otherwise, $F$ makes sense only when considered with values in $\Mm_A$ and the problem is not well posed. For convenience, we introduce the

$(*)$ condition: we say that the $(*)$ condition is fulfilled if the canonical morphism $(N\square_\Dd\Cc)\otimes\Cc\stackrel{\sim}{\rightarrow} N\square_\Dd(\Cc\otimes\Cc)$ is an isomorphism for all $N\in\Mm^\Dd$.

Recall (for example from \cite{CMZ} or \cite{McL}) that $F$ is a left adjoint to $U$ (so we have an adjointnes of functors $(F,U)$) if and only if there are natural transformations $\eta:\Id_{\Mm^\Dd}\longrightarrow UF$ and $\varepsilon:FU\longrightarrow \Id_{\Mm^\Cc}$ such that 
\begin{eqnarray}
\varepsilon_{F(N)}\circ F(\eta_N) & = & \Id_{F(N)}, \,\,\,\forall N\in\Mm^\Dd \label{1}\\
U(\varepsilon_M)\circ\eta_{U(M)} & = & \Id_{U(M)}, \,\,\,\forall M\in\Mm^\Cc  \label{2}
\end{eqnarray}
or, as $U(M)=M$ for all $M\in \Mm^\Cc$ and $U(f)=f$ for any morphism of right $\Cc$-comodules $f$, equation (\ref{2}) rewrites to 
\begin{eqnarray}
\varepsilon_M\circ\eta_{M} & = & \Id_{M}, \,\,\,\forall M\in\Mm^\Cc.  \label{3}
\end{eqnarray}
It is then natural to try to compute the (sets of) natural transformations $\Nat(\Id_{\Mm^\Dd},UF)$ and $\Nat(FU,\Id_{\Mm^\Cc})$.

\begin{proposition}\label{1.1}
Assume the ($*$) condition is satisfied. Then, with the above notations, we have $\Nat(\Id_{\Mm^\Dd},UF)\simeq {}^\Dd\Hom^\Dd(D,C)$ with the following two inverse to each other applications:
$$\Nat(\Id_{\Mm^\Dd},UF)\ni\eta\longmapsto\iota\circ\eta_\Dd\in{}^\Dd\Hom^\Dd(\Dd,\Cc)$$
$${}^\Dd\Hom^\Dd(\Dd,\Cc)\ni\alpha\longmapsto\eta(\alpha)\in\Nat(\Id_{\Mm^\Dd},UF)$$
where for $\alpha\in{}^\Dd\Hom^\Dd(\Dd,\Cc)$, $\eta(\alpha)_N(n)=n_0\otimes \alpha(n_1)$, for all $N\in\Mm^\Dd$ and $\iota:\Dd\square_{\Dd}\Cc\stackrel{\sim}{\longrightarrow} \Cc$ is the $(\Dd,\Cc)$-bicomodule isomorphism given by $\iota(d\otimes c)=\varepsilon_\Dd(d)c$ with inverse $\iota^{-1}(c)=\lambda(c_1)\otimes c_2$.
\end{proposition}
\begin{proof}
Take $\eta\in\Nat(\Id_{\Mm^\Dd},UF)$. Let $N\in \Mm^\Dd$. For each $n\in N$ define $f_n:\Dd\rightarrow N\otimes \Dd$ by $f_n(d)=n\otimes d$. Then it is easy to see that $f_n$ is a morphism of right $\Dd$-comodules, where the right $\Dd$-comodule structure on $N\otimes \Dd$ is given by $n\otimes d\mapsto n\otimes d_1\otimes d_2$. Also, for $n\in N$ denote $\eta_N(n)=n^0\otimes n^1\in N\square_{\Dd}\Cc$ where the summation symbol is again omitted. Fix $n\in N$ and suppose $\rho_N(n)=n_0\otimes n_1=\sum\limits_{i}n_i\otimes d_i\in N\otimes \Dd$. We have that $\rho_N$ is a morphism of right $\Dd$-comodules. By the naturality of $\eta$, we have the following commutative diagram:
\begin{diagram}
N 						& \rTo^{\eta_N} 									& N\square_{\Dd}\Cc \\
\dTo^{\rho_N} &																	& \dTo_{\rho_N\square_{\Dd}\Cc}\\
N\otimes \Dd  & \rTo_{\eta_{N\otimes\Dd}} & (N\otimes \Dd)\square_{\Dd}\Cc \simeq N\otimes(\Dd\square_{\Dd}\Cc)\\
\uTo^{f_{n_i}}&																	& \uTo_{f_{n_i}\square_{\Dd}\Cc}\\
\Dd						& \rTo_{\eta_\Dd}									& \Dd\square_{\Dd}\Cc\simeq \Cc
\end{diagram}
The upper diagram shows that $(\rho_N\square_{\Dd}\Cc)(n^0\otimes n^1)=\eta_{N\otimes\Dd}(n_0\otimes n_1)=\sum\limits_{i}\eta_{N\otimes\Dd}(n_i\otimes c_i)$ and the lower diagram applied for each $f_{n_i}$ and each $d_i\in \Dd$ yields $\eta_{N\otimes\Dd}(n_i\otimes d_i)=\eta_{N\otimes\Dd}(f_{n_i}(d_i))=f_{n_i}\square_{\Dd}\Cc(\eta_D(d_i))=n_i\otimes\eta_{\Dd}(d_i)$ and therefore we get 
\begin{eqnarray}
(n^0)_0\otimes (n^0)_1\otimes n^1=\sum\limits_{i}n_i\otimes \eta_\Dd(d_i). \label{4}
\end{eqnarray}
Then if we denote by $\alpha=\iota\circ\eta_\Dd$ we have $\eta_\Dd=\iota^{-1}\circ\alpha$, so (\ref{4}) becomes
\begin{eqnarray}
(n^0)_0\otimes (n^0)_1\otimes n^1=n_0\otimes \lambda(\alpha(n_1)_1)\otimes \alpha(n_1)_2. \label{5}
\end{eqnarray} 
Therefore by applying $N\otimes \varepsilon_{\Dd}\otimes \Cc$ in equation (\ref{5}) we get 
\begin{eqnarray*}
\eta_N(n)=n^0\otimes n^1 & = & (n^0)_0\varepsilon_\Dd((n^0)_1)\otimes n^1 \\
& = & n_0\varepsilon_\Dd(\lambda(\alpha(n_1)_1)\otimes \alpha(n_1)_2\\
& = & n_0\varepsilon_\Cc(\alpha(n_1)_1)\otimes \alpha(n_2)_2\\
& = & n_0\otimes \varepsilon_\Cc(\alpha(n_1)_1)\alpha(n_1)_2\\
& = & n_0\otimes\alpha(n_1) \,\,\,({\rm by\,the\,counit\,property})\\
\end{eqnarray*}
As $\iota$ is a morphism of right $\Cc$-comodules, it is also a morphism of right $\Dd$-comodules, and therefore $\alpha=\iota\circ\eta_D$ is also a morphism in $\Mm^\Dd$. Now the relation $\eta_D=\iota^{-1}\circ\alpha$ rewrites 
\begin{eqnarray}
d_1\otimes \alpha(d_2)=\eta_D(d) & = &\iota^{-1}(\alpha(d))=\lambda(\alpha(d)_1)\otimes\alpha(d)_2 \label{6}
\end{eqnarray}
and this is exactly the fact that $\alpha$ is a morphism of left $\Dd$-comodules. Therefore $\alpha\in{}^\Dd\Hom^\Dd(\Dd,\Cc)$. 

Conversely, start with such an $\alpha$. Then relation (\ref{6}) holds and consequently for $N\in\Mm^\Dd$ and $n\in N$, we have $n_{00}\otimes n_{01}\otimes \alpha(n_1)=n_0\otimes n_1\otimes \alpha(n_2)=n_0\otimes\lambda(\alpha(n_1)_1)\otimes \alpha(n_1)_2$ (by applying the coassociativity of the comultiplication and (\ref{6}) for $n_1$, the second position of $\rho_N(n)=n_0\otimes n_1\in N\otimes\Dd$). This shows that actually $n_0\otimes \alpha(n_1)\in \Ker(\omega_{N,\Cc})$ and therefore it makes sense to define $\eta_N:N\longrightarrow N\square_{\Dd}\Cc$, $\eta_{N}(n)=n_0\otimes \alpha(n_1)$. Because $\alpha$ is a morphism of right $\Dd$-comodules we get
\begin{eqnarray*}
\rho_{N\square_{\Dd}\Cc}(\eta_N(n))& = & \eta_N(n)_0\otimes\eta_N(n)_1=\rho_{N\square_{\Dd}\Cc}(n_0\otimes\alpha(n_1))\\
& = & n_0\otimes \alpha(n_1)_1\otimes\lambda(\alpha(n_1)_2) \,\,\,({\rm the\,right\,\Dd\,comodule\,structure\,on\,N\square_{\Dd}\Cc}\\
&   & {\rm\,comes\,from\,the\,right\,\Cc\,comodule\,structure\,via\,\lambda:\Cc\longrightarrow \Dd})\\
& = & n_0\otimes\alpha(n_1)\otimes n_2 \,\,\,({\rm as\,\alpha\,is\,a\,morphism\,in\,\Mm^\Dd})\\
& = & \eta_N(n_0)\otimes n_1
\end{eqnarray*}
This shows that $\eta_N$ is a morphism of right $\Dd$-comodules. Now note that $\eta$ is a natural transformation: indeed if $f:N\longrightarrow N'$ is a morphism of right $\Dd$-comodules we need to show that the following diagram is commutative:
\begin{diagram}
N        & \rTo^{\eta_N}    & N\square_{\Dd}\Cc \\
\dTo^{f} &                  & \dTo_{f\square_{\Dd}\Cc}    \\
N'       & \rTo_{\eta_{N'}} & N'\square_{\Dd}\Cc
\end{diagram}
This follows as for $n\in N$
\begin{eqnarray*}
\eta_{N'}(f(n)) & = & f(n)_0\otimes \alpha(f(n)_1)\\
& = & f(n_0)\otimes \alpha(n_1)\,\,\,({\rm because\,}f(n)_0\otimes f(n)_1=f(n_0)\otimes n_1\,{\rm as\,}\\
&   & f\,{\rm is\,a\,morphism\,in\,}\Mm^\Dd)\\
& = & (f\square_{\Dd}\Cc)(\eta_N(n))
\end{eqnarray*}
Now note that for $\eta\in\Nat(\Id_{\Mm^\Dd},UF)$ we have $\eta_N(n)=n_0\otimes\alpha(n_1)$ and therefore $\eta=\eta(\alpha)$ with the notations in the statement of the proposition. Thus 
\begin{eqnarray}
\eta=\eta(\iota\circ\eta_D)\label{7}
\end{eqnarray} 
Also for $\alpha\in{}^\Dd\Hom^\Dd(\Dd,\Cc)$, $\eta(\alpha)_N(n)=n_0\otimes\alpha(n_1)$ and then we get that $(\iota\circ\eta_{D}(\alpha))(d)=\iota(d_1\otimes\alpha(d_2))=\varepsilon_\Dd(d_1)\alpha(d_2)=\alpha(\varepsilon_\Dd(d_1)d_2)=\alpha(d)$. Therefore
\begin{eqnarray}
\iota\circ\eta_{D}(\alpha)=\alpha \label{8}
\end{eqnarray}
Now equations (\ref{7}) and (\ref{8}) show that the applications in the statement of the proposition are inverse to each other, and the proof is finished.
\end{proof}

For each natural transformation $\epsilon\in\Nat(FU,\Id_{\Mm^\Cc})$ we can associate $\beta\in\Hom^\Cc(\Cc\square_{\Dd}\Cc,\Cc)$ by putting $\beta=\epsilon_\Cc$. Under certain conditions this also becomes a morphism of left $C$-comodules. In fact, for $\Cc\square_{\Dd}\Cc$ to have a left $\Cc$-comodule structure we need an isomorphism of $A$-(bi)modules $\Cc\otimes(\Cc\square_{\Dd}\Cc)\simeq (\Cc\otimes\Cc)\square_{\Dd}\Cc$. Therefore in the next Proposition we will assume that $\omega=\omega_{\Cc,\Cc}$ is pure in ${}_A\Mm$. This is not a very restrictive condition as it will be seen that it holds many situations including in all the cases of our applications. Moreover, this condition is automatically fulfilled when the functors $(F,U)$ form a Frobenius pair. If this does not hold, it will be seen from the proof of the following proposition that it is difficult to describe these natural transformations.

\begin{proposition}\label{1.2}
With the above notations, if $\omega_{\Cc,\Cc}$ is pure in ${}_A\Mm$ then $\Nat(FU,\Id_{\Mm^\Cc})\simeq {}^\Cc\Hom^\Cc(\Cc\square_{\Dd}\Cc,\Cc)$, where the applications giving the equivalence are given by
$$\Nat(FU,\Id_{\Mm^\Cc})\ni\epsilon\longmapsto \epsilon_{\Cc}\in{}^\Cc\Hom^\Cc(\Cc\square_{\Dd}\Cc,\Cc)$$
and
$${}^\Cc\Hom^\Cc(\Cc\square_{\Dd}\Cc,\Cc)\ni\beta\longmapsto\epsilon(\beta)\in\Nat(FU,\Id_{\Mm^\Cc})$$
where for $\beta\in{}^\Cc\Hom^\Cc(\Cc\square_{\Dd}\Cc,\Cc)$, $\epsilon(\beta)_M(m\otimes c)=m_0\varepsilon_\Cc(\beta(m_1\otimes c))$ for all $M\in\Mm^\Cc$ and $m\otimes c\in M\square_{\Dd}\Cc$ (here $m_0\varepsilon_\Cc(\beta(m_1\otimes c))$ means $M\otimes \varepsilon_\Cc\beta$ applied to $m_0\otimes m_1\otimes c\in (M\otimes C)\square_\Dd\Cc\simeq M\otimes (\Cc\square_\Dd\Cc)$ composed to the canonical isomorphism $M\otimes A\simeq A$).
\end{proposition}
\begin{proof}
Let $M$ be a right $\Cc$-comodule and $\epsilon\in\Nat(FU,\Id_{\Mm^\Cc})$. Denote $\beta=\epsilon_\Cc$. Because $\omega_{\Cc,\Cc}$ is pure in ${}_A\Mm$, we always have a natural isomorphism of right $A$-modules $\psi_M:M\otimes(\Cc\square_{\Dd}\Cc)\stackrel{\sim}{\longrightarrow}(M\otimes \Cc)\square_{\Dd}\Cc$ (see \cite{BW}, Section 21, 21.4). Pick $m\otimes c\in M\square_{\Dd}\Cc$. By the above mentioned isomorphism $\psi_M$, there are $m_i\in M$, $c_i\otimes c'_i\in\Cc\square_{\Dd}\Cc$ such that $\psi_M(\sum\limits_im_i\otimes c_i\otimes c'_i)=m_0\otimes m_1\otimes c\in(M\otimes \Cc)\square_{\Dd}\Cc$. We convey to identify each element from $M\otimes(\Cc\square_{\Dd}\Cc)$ with its image via $\psi_M$, for convenience. Then by the naturality of $\epsilon$ we have a commutative diagram
\begin{diagram}
																		& M\square_{\Dd}\Cc               & \rTo^{\epsilon_M}							 & M \\
																		& \dTo^{\rho_M\square_{\Dd}\Cc}   &                  							 & \dTo_{\rho_M} \\
M\otimes(\Cc\square_{\Dd}\Cc)\simeq & (M\otimes \Cc)\square_{\Dd}\Cc  & \rTo_{\epsilon_{M\otimes \Cc}} & M\otimes \Cc  \\
																		& \uTo^{g_{m_i}\square_{\Dd}\Cc}  &                                & \uTo_{g_{m_i}}\\
																		& \Cc\square_{\Dd}\Cc   					& \rTo_{\beta}									 & \Cc
\end{diagram}
where for each $i$, $g_{m_i}:\Cc\longrightarrow M\otimes\Cc$ is the right $\Cc$-comodule morphism defined by $g_{m_i}(c)=m_i\otimes c$. By the commutativity of the lower part of the diagram, for each $i$ we get $\epsilon_{M\otimes\Cc}(m_i\otimes c_i\otimes c'_i)=\epsilon_{M\otimes \Cc}(g_{m_i}\square_\Dd\Cc(c_i\otimes c'_i))=g_{m_i}(\beta(c_i\otimes c'_i))=m_i\otimes\beta(c_i\otimes c'_i)$, so then 
\begin{eqnarray}
\epsilon_{M\otimes \Cc}(\sum\limits_im_i\otimes c_i\otimes c'_i) & = &\sum\limits_im_i\otimes\beta(c_i\otimes c'_i) \label{9}
\end{eqnarray} 
By the upper part of the diagram we have $\epsilon_{M \otimes\Cc}(m_0\otimes m_1\otimes c)=\epsilon_{M \otimes\Cc}(\rho_M\square_{\Dd}\Cc(m\otimes c))=\rho_M(\epsilon_M(m\otimes c))$ and therefore 
\begin{eqnarray}
\epsilon_{M\otimes\Cc}(m_0\otimes m_1\otimes c) & =& \epsilon_M(m\otimes c)_0\otimes \epsilon_M(m\otimes c)_1 \label{10}
\end{eqnarray}
Combining (\ref{9}) and (\ref{10}) and keeping in mind the identification between $M\otimes(\Cc\square_{\Dd}\Cc)$ and $(M\otimes \Cc)\square_{\Dd}\Cc$ made via $\psi_M$ we get
\begin{eqnarray}
m_0\otimes\beta(m_1\otimes c) & = & \epsilon_M(m\otimes c)_0\otimes \epsilon_M(m\otimes c)_1 \label{11}
\end{eqnarray}
and therefore by applying $\varepsilon_\Cc$ on the second position we get 
\begin{eqnarray}
\epsilon_M(m\otimes c) & = & m_0\varepsilon_\Cc(\beta(m_1\otimes c)) \label{12}
\end{eqnarray}
(this formula is aways understood as $\epsilon_M(m\otimes c)\,=\,\sum\limits_im_i\varepsilon_\Cc(\beta(c_i\otimes c'_i))$, where $\sum\limits_im_i\otimes c_i\otimes c'_i\psi_M^{-1}(m_0\otimes m_1\otimes c)$).\\
We have that $\beta=\epsilon_\Cc$ is a morphism of right $\Cc$-comodules. Now writing equation (\ref{11}) for $M=\Cc$ and $m\otimes c=c'\otimes c\in \Cc\square_{\Dd}\Cc$ we get $c'_1\otimes \beta(c'_2\otimes c)=\beta(c'\otimes c)_1\otimes\beta(c'\otimes c)_2$ and this shows that $\beta$ is a morphism in ${}^\Cc\Mm$ (we have already seen that $\Cc\square_{\Dd}\Cc$ has a left $\Cc$-comodule structure because $\omega_{\Cc,\Cc}$ is pure in ${}_A\Mm$). \\
Conversely, take $\beta\in{}^\Cc\Hom^\Cc(\Cc\square_{\Dd}\Cc,\Cc)$. For $M\in\Mm^\Cc$, as $\psi_M$ is an isomorphism for all $M$ we can consider the application $\epsilon_M:M\square_{\Dd}\Cc\longrightarrow M$ given by the formula in (\ref{12}), so $\epsilon_M=\nu_M\circ(M\otimes \varepsilon_\Cc)\circ(M\otimes \beta)\circ\psi_M^{-1}\circ(\rho_M\square_{\Dd}\Cc)$, where $\nu_M:M\otimes A\stackrel{\sim}{\rightarrow}M$ is the canonical isomorphism. Then $\epsilon_M$ is a morphism of right $\Cc$-comodules: 
\begin{eqnarray*}
\epsilon_M(m\otimes c)_0\otimes\epsilon_M(m\otimes c)_1 & = & m_{00}\otimes m_{01}\varepsilon_\Cc(\beta(m_1\otimes c))\\
& = & m_0 \otimes m_{11}\varepsilon_\Cc(\beta(m_{12}\otimes c))\;\;\;({\rm the\,coassociativity\,of\,}\rho_M)\\
& = & m_0\otimes \beta(m_1\otimes c)_1\varepsilon_\Cc(\beta(m_1\otimes c)_2) \;\;\;({\rm because\,}\beta\\
&   & {\rm is\,a\,morphism\,in\,}{}^\Cc\Mm)\\
& = & m_0\otimes \beta(m_1\otimes c)\\
& = & m_0\otimes \varepsilon_\Cc(\beta(m_1\otimes c)_1)\beta(m_1\otimes c)_2\\
& = & m_0\varepsilon_\Cc(\beta(m_1\otimes c)_1)\otimes \beta(m_1\otimes c)_2\\
& = & m_0\varepsilon_\Cc(\beta(m_1\otimes c_1))\otimes c_2 \;\;\;({\rm as\,}\beta\,{\rm is\,a\,morphism\,in\,}\Mm^\Cc)\\
& = & \epsilon_M(m\otimes c_1)\otimes c_2
\end{eqnarray*}
Also, note that $\epsilon$ is a natural transformation, that is, for a morphism of right $\Cc$-comodules $f:M\longrightarrow M'$ the following diagram is commutative:
\begin{diagram}
M\square_{\Dd}\Cc 			 & \rTo^{\epsilon_M}    & M\\
\dTo^{f\square_{\Dd}\Cc} &									    & \dTo_{f} \\
M'\square_{\Dd}\Cc       & \rTo_{\epsilon_{M'}} & M'
\end{diagram}
Indeed for $m\otimes c\in M\square_{\Dd}\Cc$ we have
\begin{eqnarray*}
\epsilon_{M'}(f\square_{\Dd}\Cc(m\otimes c)) & = & \epsilon_{M'}(f(m)\otimes c)\\
& = & f(m)_0\varepsilon_\Cc(\beta(f(m)_1\otimes c)) \\
& = & f(m_0)\varepsilon_\Cc(\beta(m_1\otimes c)) \;\;\;({\rm because\,}f\,{\rm is\,a\,morphism\,in\,}\Mm^\Cc)\\
& = & f(m_0\varepsilon_\Cc(\beta(m_1\otimes c)))\\
& = & f(\epsilon_M(m\otimes c))
\end{eqnarray*}
Thus $\epsilon=\epsilon(\beta)\in\Nat(FU,\Id_{\Mm^\Cc})$ is a well defined natural transformation. Now, if $\epsilon\in\Nat(FU,\Id_{\Mm^\Cc})$ then by formula (\ref{12}) we see that $\epsilon=\epsilon(\beta)$ for $\beta=\epsilon_\Cc$ and therefore 
\begin{equation}
\epsilon=\epsilon(\epsilon_\Cc) \label{13}
\end{equation}
Also, if $\beta\in{}^\Cc\Hom^\Cc(\Cc\square_{\Dd}\Cc,\Cc)$, then for $c\otimes c'\in\Cc\square_{\Dd}\Cc$ we have
\begin{eqnarray*}
\epsilon(\beta)_\Cc(c\otimes c') & = & c_1\varepsilon_\Cc(\beta(c_2\otimes c'))\\
& = & \beta(c\otimes c')_1\varepsilon_{\Cc}(\beta(c\otimes c'))\;\;\;({\rm because\,}\beta\,{\rm is\,a\,morphism\,in\,}\Mm^\Cc)\\
& = & \beta(c\otimes c')
\end{eqnarray*}
and therefore 
\begin{equation}
\epsilon(\beta)_\Cc=\beta \label{14}
\end{equation}
Equations (\ref{13}) and (\ref{14}) show that the applications given in the statement of the Proposition are inverse to each other and the proof is finished.
\end{proof}

\begin{definition}
Let $\lambda:\Cc\longrightarrow \Dd$ be a morphism of corings. We say $\Cc\stackrel{\lambda}{\longrightarrow}\Dd$ is a right Frobenius extension of corings (or $\lambda$ is Frobenius) if the functor $U:\Mm^{\Cc}\longrightarrow \Mm^\Dd$ is a Frobenius functor. 
\end{definition}

In the case the categories $\Mm^\Cc$ and $\Mm^\Dd$ are abelian categories, then if $(F,U)$ form an adjoint pair (i.e. $F$ is a left adjoint for $U$) then $F$ is right exact. We can see that this is also true for more general categories, in particular for the case we are interested in. Recall from \cite{BW} that the category of right $\Cc$-comodules $\Mm^\Cc$ has cokernels and the image of every morphism of right $\Cc$-comodules $f:M\longrightarrow M'$ is a subcomodule in $M'$. We say that a sequence $M'\stackrel{u}{\longrightarrow}M\stackrel{v}{\longrightarrow}M''{\longrightarrow}0$ is exact provided that $v$ is an epimorphism, $v$ has Kernel and ${\rm Im}(u)=\Ker(v)$. A certain right exactness of $F$ that we need can be proved here without the assumption that $\Mm^\Cc$ and $\Mm^\Dd$ are abelian, which in turn is used to prove that the condition $\omega_{\Cc,\Cc}$ is pure in ${}_A\Mm$ holds when $\Cc\stackrel{\lambda}{\longrightarrow}\Dd$ is a Frobenius extension. Following the ideas from \cite{BW} 21.5, 18.16, 3.19 and 3.3, we can prove the following

\begin{proposition}\label{1.4}
Suppose the $(*)$-condition holds and that $(F,U)$ form an adjoint pair, i.e., $F$ is a left adjoint to $U$. Then $\omega_{\Cc,\Cc}$ is pure in ${}_A{\Mm}$.
\end{proposition}
\begin{proof}
By \cite{BW}, Section 21 (21.4, Tensor-cotensor relations) we see that $\omega_{\Cc,\Cc}$ is pure in ${}_A\Mm$ if and only if the canonical morphism $\psi_M:M\otimes(\Cc\square_{\Dd}\Cc)\longrightarrow (M\otimes \Cc)\square_{\Dd}\Cc$ is an isomorphism for any $M\in{}_A\Mm$. Let $F_2\longrightarrow F_1\longrightarrow M\longrightarrow 0$ be an exact sequence of $A$-modules with free $A$-modules $F_2$ and $F_1$. Then we have the commutative diagram
\begin{diagram}
(F_2\otimes \Cc)\square_{\Dd}\Cc & \rTo & (F_1\otimes \Cc)\square_{\Dd}\Cc & \rTo & (M\otimes \Cc)\square_{\Dd}\Cc & \rTo & 0\\
\uTo_{\simeq}^{\psi_{F_2}} & & \uTo_{\simeq}^{\psi_{F_1}} & & \uTo & & \\
F_2\otimes(\Cc\square_{\Dd}\Cc) & \rTo & F_1\otimes(\Cc\square_{\Dd}\Cc) & \rTo & M\otimes(\Cc\square_{\Dd}\Cc) & \rTo 0
\end{diagram}
where the bottom row is exact by the right exactness of the tensor product and the first two vertical arrows are isomorphisms because the tensor and cotensor product commute with coproducts. Then in order to finish the proof, by the above diagram it is enough to prove that the top row is exact as a sequence of $A$-modules, that is, the sequence $F(F_2\otimes \Cc)\longrightarrow F(F_1\otimes \Cc)\longrightarrow F(M\otimes \Cc)\longrightarrow 0$ is an exact sequence of $A$-modules ($M\otimes \Cc=U(M\otimes \Cc)$ is considered here as a $\Dd$-comodule). Denote $Z=F_2\otimes \Cc$, $Y=F_1\otimes \Cc$, $X=M\otimes\Cc$. Because $F$ is a left adjoint to $U$, for every right $\Cc$-comodule $W$ we have a commutative diagram
\begin{diagram}
& 0 & \rTo & \Hom^\Cc(F(X),W) & \rTo & \Hom^\Cc(F(Y),W) & \rTo & \Hom^\Cc(F(Z),W)\\
(D1)& & & \dTo_{\simeq} & & \dTo_{\simeq} & & \dTo_{\simeq} \\
& 0 & \rTo & \Hom^\Dd(X,U(W)) & \rTo & \Hom^\Dd(Y,U(W)) & \rTo & \Hom^\Dd(Z,U(W)) 
\end{diagram}
We prove that the bottom row of this last diagram is commutative for every $W\in\Mm^\Cc$, and this will imply that the top row is exact too. By \cite{BW}, Section 3 (3.3) we have that for any two right $\Dd$-comodules $N$ and $N'$ there is an exact sequence of abelian groups
$$0\longrightarrow \Hom^\Dd(N,N')\longrightarrow \Hom_A(N,N')\stackrel{\nu_{N,N'}}{\longrightarrow} \Hom_A(N,N'\otimes \Dd)$$
where $\nu_{N,N'}(f)=\rho_{N'}\circ f-(f\otimes \Cc)\circ \rho_N$. Therefore, as in \cite{BW}, 3.19, we have a commutative diagram yielded by the exactness of the sequence of $A$ modules 
$Z\longrightarrow Y\longrightarrow X\longrightarrow 0$ which is easy to see that it is exact also as a sequence of right $\Dd$-comodules (because it is exact in $\Mm_A$; note that we do not need to have that ${}_A\Dd$ is flat so that the category ${}_A\Mm$ to be abelian).
\begin{diagram}
& & 0 & & 0 & & 0\\
& & \dTo & & \dTo & & \dTo\\
0 & \rTo & \Hom^\Dd(X,U(W)) & \rTo & \Hom^{\Dd}(Y,U(W)) & \rTo & \Hom^\Dd(Z,U(W)) &\\ 
& & \dTo & & \dTo & & \dTo &\\
0 & \rTo & \Hom_A(X,U(W)) & \rTo & \Hom_A(Y,U(W)) & \rTo & \Hom_A(Z,U(W)) &\\
& & \dTo & & \dTo & & \dTo &\\
0 & \rTo & \Hom_A(X,U(W)\otimes\Dd) & \rTo & \Hom_A(Y,U(W)\otimes\Dd) & \rTo & \Hom_A(Z,U(W)\otimes\Dd) &
\end{diagram}

The two bottom rows are exact by the properties of the $\Hom_A(-,-)$ functor and the vertical columns are exact by the previous observation (\cite{BW}, 3.3). Then it follows that the top row of this last diagram exact (in ${}_A\Mm$), that is, the bottom row of (D1) is exact and therefore the first row in diagram (D1) is exact too. Now denote $F(Z)\stackrel{v}{\longrightarrow}F(Y)\stackrel{u}{\longrightarrow}F(X)$ the morphisms induced by $F$ from the exact sequence $Z\longrightarrow Y\longrightarrow X\longrightarrow 0$. Put $W=F(X)/{\rm Im}(u)$ (which has a natural $\Cc$-comodule structure as ${\rm Im}(u)$ is a subcomodule of $F(X)$!) and $\pi:F(X)\longrightarrow F(X)/{\rm Im}(u)$ the canonical projection. As the first row in diagram (D1) is exact, we see that $\Hom^\Cc(u,W)$ is injective. But $\Hom^\Cc(u,W)(\pi)=\pi\circ u=0$, so $\pi=0$ and therefore ${\rm Im}(u)=F(X)$. Now as ${\rm Im}(v)$ is a subcomodule of $F(Y)$, for the right $\Cc$-comodule $W=F(Y)/{\rm Im}(v)$ and the canonical projection $p:F(Y)\longrightarrow F(Y)/{\rm Im}(v)$, from the exactness in (D1) we again find that $\Ker(\Hom^\Cc(v,W))={\rm Im}(\Hom^\Cc(u,W))$. Then $\Hom^\Cc(v,W)(p)=p\circ v=0$, so $p\in\Ker(\Hom^\Cc(v,W))$ and therefore there is $h:F(X)\longrightarrow W$ such that $h\circ u=p$. Hence if $y\in F(Y)$ is such that $u(y)=0$, we get $h(u(y))=0$ so $p(y)=0$, i.e. $y\in {\rm Im}(v)$ showing that ${\rm Im}(v)\subseteq\Ker(u)$. The converse inclusion is obvious as $u\circ v=0$ (it follows from the functoriality of $F$ and the exactness of the sequence $Z\longrightarrow Y\longrightarrow X \longrightarrow 0$). With this, we get that the sequence 
$$(F_2\otimes \Cc)\square_{\Dd}\Cc \longrightarrow(F_1\otimes \Cc)\square_{\Dd}\Cc \longrightarrow (M\otimes \Cc)\square_{\Dd}\Cc \longrightarrow 0$$
is exact (in ${}_A\Mm$) and the proof is finished as shown before.
\end{proof}

\begin{theorem}\label{th.Frob}
Suppose the ($*$)-condition holds for the coring morphism $\lambda$. Then $\Cc\stackrel{\lambda}{\longrightarrow}\Dd$ is a right Frobenius extension of corings if and only if $\omega_{\Cc,\Cc}$ is pure in ${}_A\Mm$ and there are $\alpha\in{}^\Dd\Hom^\Dd(\Dd,\Cc)$ and $\beta\in{}^\Cc\Hom^\Cc(\Cc\square_{\Dd}\Cc,\Cc)$ such that 
\begin{equation}\label{15}
\beta(c_1\otimes\alpha\lambda(c_2)) \, = \, \beta(\alpha\lambda(c_1)\otimes c_2) \, = c
\end{equation} 
for all $c\in \Cc$.
\end{theorem}
\begin{proof}
Assume $\Cc\stackrel{\lambda}{\longrightarrow}\Dd$ is a Frobenius extension. Then by Proposition \ref{1.4} we have that $\omega_{\Cc,\Cc}$ is pure in ${}_A\Mm$. As $(F,U)$ is a Frobenius pair of adjoint functors, we have that $F$ is also a left adjoint to $U$. Let $\eta:\Id_{\Mm^\Dd}\longrightarrow UF$ and $\epsilon:FU\longrightarrow \Id_{\Mm^\Cc}$ be the unit and counit of this adjunction, thus satisfying the equations (\ref{1}) and (\ref{2}). Then note that Propositions \ref{1.1} and \ref{1.2} apply, and then we obtain $\alpha\in{}^\Dd\Hom^\Dd(\Dd,\Cc)$ and $\beta\in{}^\Cc\Hom^\Cc(\Cc\square_{\Dd}\Cc,\Cc)$ such that for $N\in\Mm^\Dd$ and $n\in N$ 
\begin{equation}\label{16}
\eta_N(n)=n_0\otimes\alpha(n_1)
\end{equation}
and for $M\in\Mm^\Cc$ and $m\otimes c\in M\square_\Dd \Cc$ 
\begin{equation}\label{17}
\epsilon_M(m\otimes c)=m_0\varepsilon_\Cc(\beta(m_1\otimes c))
\end{equation}
By the naturality of $\epsilon$ we have 
\begin{equation}\label{18}
\beta\circ(\iota\square_{\Dd}\Cc)=\iota\circ\epsilon_{F(\Dd)}
\end{equation}
where $\iota:\Dd\square_{\Dd}\Cc\stackrel{\sim}{\longrightarrow}\Cc$ is the isomorphism from Proposition \ref{1.1}. Then for $c\in\Cc$ we have: 
\begin{eqnarray*}
c & = & \iota(\lambda(c_1)\otimes c_2) = \iota(\epsilon_{F(\Dd)}\circ F(\eta_\Dd)(\lambda(c_1)\otimes c_2)) \;\;\;({\rm by\,equation\,}(\ref{1}))\\
  & = & \iota(\epsilon_{F(\Dd)}(\lambda(c_1)_1\otimes \alpha(\lambda(c_1)_2)\otimes c_2)) \;\;\;({\rm by\,}(\ref{16}))\\
  & = & \beta(\iota\square_{\Dd}\Cc(\lambda(c_1)_1\otimes \alpha(\lambda(c_1)_2)\otimes c_2))\;\;\;({\rm from\,}(\ref{18}))\\
  & = & \beta(\iota\square_{\Dd}\Cc(\lambda(\alpha(\lambda(c_1))_1)\otimes\alpha(\lambda(c_1))_2\otimes c_2)) \;\;\;({\rm as\,}\alpha\,{\rm is\,a\,morphism\,in\,}{}^\Dd\Mm)\\
  & = & \beta(\varepsilon_\Dd(\lambda(\alpha\lambda(c_1)_1))\alpha\lambda(c_1)_2\otimes c_2) \\
  & = & \beta(\varepsilon_\Cc(\alpha\lambda(c_1)_1)\alpha\lambda(c_1)_2\otimes c_2) \\
  & = & \beta(\alpha\lambda(c_1)\otimes c_2)
\end{eqnarray*}
Also 
\begin{eqnarray*}
c & = & \epsilon_\Cc(\eta_\Cc(c))\;\;\;({\rm by\,}(\ref{3})\,{\rm for\,}M=\Cc) \\
  & = & \beta(\eta_\Cc(c)) \\
  & = & \beta(c_1\otimes\alpha\lambda(c_2)) \;\;\;({\rm by\,}(\ref{16}))
\end{eqnarray*}
and therefore (\ref{15}) is proved. 
Conversely, assume there are $\alpha$ and $\beta$ such that (\ref{15}) holds and also that $\omega_{\Cc,\Cc}$ is pure in ${}_A\Mm$.  Then again by Propositions \ref{1.1} and \ref{1.2} we can find the natural transformations $\eta:\Id_{\Mm^\Dd}\longrightarrow UF$ and $\epsilon:FU\longrightarrow \Id_{\Mm^\Cc}$ such that conditions (\ref{16}) and (\ref{17}) are fulfilled. Then, with notations as above and $n\otimes c\in N\square_{\Dd}\Cc$ we have
\begin{eqnarray*}
\epsilon_{F(N)}(F(\eta_N)(n\otimes c)) & = & \epsilon_{F(N)}(n_0\otimes \alpha(n_1)\otimes c) \;\;\;{\rm by\,}(\ref{16})\\
 & = & n_0\otimes\alpha(n_1)_1\varepsilon_\Cc(\beta(\alpha(n_1)_2\otimes c))    \;\;\;({\rm from\,}(\ref{17}){\,\rm for\,}M=F(N))\\
 & = & n_0\otimes\beta(\alpha(n_1)\otimes c)_1\varepsilon_\Cc(\beta(\alpha(n_1)\otimes c)_2) \;\;\;({\rm because\,}\beta\in{}^\Cc\Mm)\\
 & = & n_0\otimes\beta(\alpha(n_1)\otimes c) \\
 & = & n\otimes\beta(\alpha\lambda(c_1)\otimes c_2) \;\;\;{\rm as\,}n\otimes c\in N\square_{\Dd}\Cc\\
 & = & n\otimes c \;\;\;{\rm by\,}(\ref{15})\\
\end{eqnarray*}
It is not difficult to see that the compositions of functions involved above make sense and therefore these computations yield (\ref{1}). Also for $M\in\Mm^\Cc$ and $m\in M$ we have 
\begin{eqnarray*}
\epsilon_M(\eta_M(m)) & = & \epsilon_M(m_0\otimes\alpha\lambda(m_1)) \;\;\;({\rm the\,right\,}\Dd{\rm-comodule\,structure\,of\,}M)\\
 & = & m_{00}\varepsilon_\Cc\beta(m_{01}\otimes \alpha\lambda(m_1)) \\
 & = & m_{0}\varepsilon_\Cc\beta(m_{11}\otimes \alpha\lambda(m_{12})) \\
 & = & m_0\varepsilon_\Cc(m_1) \;\;\;{\rm by\,}(\ref{15}) \\
 & = & m
\end{eqnarray*}
and then (\ref{2}) holds. Equations (\ref{1}) and (\ref{2}) show that $F$ is also a left adjoint to $U$ which amounts to the fact that $F$ is a right adjoint to $U$, showing that $\Cc\stackrel{\lambda}{\longrightarrow}\Dd$ is a Frobenius extension.
\end{proof}

\begin{remark}
If $M$ is a $\Cc$-$\Cc$ bicomodule, then we have ${}^\Cc\Hom^\Cc(M,\Cc)\simeq \{u\in{}_A\Hom_A(M,A)\mid u(m_0)m_1=m_{-1}u(m_0)\}$. Indeed, if $h\in{}^\Cc\Hom^\Cc(M,\Cc)$, then $u=\varepsilon_\Cc\circ h\in{}_A\Hom_A(M,\Cc)$ and $u(m_0)m_1=\varepsilon_\Cc(h(m_0))m_1=\varepsilon(h(m)_1)h(m)_2$ (because $h$ is a morphism of right $\Cc$-comodules) and therefore $u(m_0)m_1=h(m)$ and similarly $m_{-1}u(m_0)=h(m)$ so 
\begin{equation}\label{19}
m_{-1}u(m_0)=u(m_0)m_1 
\end{equation}
Conversely, for $u\in{}_A\Hom_A(M,A)$ which satisfies (\ref{19}), we can easily see that $h_u=(m\mapsto u(m_0)m_1=m_{-1}u(m_0))$ is a morphism of left and also of right $\Cc$-comodules because of (\ref{19}) and that $\varepsilon_\Cc\circ h=u$ if and only if $h=h_u$. Therefore we can equivalently express the statement of Theorem \ref{th.Frob}:
\begin{corollary}\label{1.cor}
If the ($*$)-condition holds for the coring morphism $\lambda$, then $\Cc\stackrel{\lambda}{\longrightarrow}\Dd$ is a right Frobenius extension of corings if and only if $\omega_{\Cc,\Cc}$ is pure in ${}_A\Mm$ and there are $\alpha\in{}^\Dd\Hom^\Dd(\Dd,\Cc)$ and $\gamma\in{}_A\Hom_A(\Cc\square_{\Dd}\Cc,A)$ such that 
\begin{eqnarray}
\gamma(\alpha\lambda(c_1)\otimes c_2) & = & \gamma(c_1\otimes\alpha\lambda(c_1))=\varepsilon_\Cc(c) \label{20}
\end{eqnarray}
and equation (\ref{19}) holds for $M=\Cc\square_{\Dd}\Cc$, that is, $c_1\gamma(c_2\otimes c')=\gamma(c\otimes c'_1)c'_2=c\otimes c'$ for $c\otimes c'\in\Cc\square_\Dd\Cc$.
\end{corollary}
\end{remark}

\section{Examples}
We consider here several situations where the above considerations apply. We also give some finiteness theorems for the cases we study.

\begin{example}\label{e2.1}
Let $I$ be a set, $(\Dd,\Delta_\Dd,\varepsilon_\Dd)$ be an $A$-coring. Let $\Cc=\bigoplus\limits_{i\in I}\Dd$ (the direct sum of $A-A$-bimodules) and denote by $\sigma_i:\Dd\longrightarrow \Cc$ be the canonical injections and $p_i:\Cc\longrightarrow \Dd$ be the canonical projection on the $i$-th component. Then any element in $\Cc$ is of the form $\sum\limits_{i\in I}\sigma_i(d_i)$ with $d_i\in\Dd$ and it is not difficult to see that $\Cc$ becomes an $A$-coring by comultiplication $\Delta_{\Cc}$ given by
$$\Delta_{\Cc}(\sum\limits_{i\in I}\sigma_i(d_i))=\sum\limits_{i\in I}\sigma_i(d_{i1})\otimes\sigma_i(d_{i2})$$
and counit $\varepsilon_\Cc$ given by the formula
$$\varepsilon_\Cc(\sum\limits_{i\in I}\sigma_i(d_i))=\sum\limits_{i\in I}\varepsilon_\Dd(d_i)$$
Let $\lambda:\Cc\longrightarrow \Dd$ be the $A$-bimodule morphism defined by $\lambda(\sum\limits_{i\in I}\sigma_i(d_i))=\sum\limits_{i\in I}d_i$. By the above definitions it is clear that $\lambda$ is a morphism of corings. By \cite{I1}, Proposition 2.6 we have that the category $\Mm^\Cc$ is equivalent to the product of categories $(\Mm^\Dd)^I$. Recall that any $\Cc$-comodule $M$ is given by $M=\bigoplus\limits_{i\in I}M_i$ where each $M_i$ is a $\Dd$-comodule and if $\theta_i:M_i\longrightarrow M$ is the canonic injection, then the $\Cc$-comodule structure of $M$ is given by $\rho_M:M\longrightarrow M\otimes\Cc$ with $\rho(\sum\limits_{i\in I}\theta_i(m_i))=\sum\limits_{i\in i}\theta_i(m_{i0})\otimes \sigma_i(m_{i1})$. The right $\Dd$-comodule structure of $M$ is given by $\rho'_M:M\longrightarrow M\otimes\Dd$ with $\rho'_M(\sum\limits_{i\in I}(\theta_i(m_i))=\sum\limits_{i\in I}\theta_i(m_{i0})\otimes m_{i1}\in M\otimes \Dd$. It is easy to see that the $\Dd$-comodule structure of $M$ as the direct sum of the $\Dd$-comodules $(M_i)_{i\in I}$ is the same as the one obtained from the right $\Cc$-comodule structure of $M$ via $\lambda$. Therefore the functor $U$ associated to the extension of corings $\Cc\stackrel{\lambda}{\longrightarrow}\Dd$ coincides to the coproduct functor from $(\Mm^\Dd)^I$ to $\Mm^\Dd$.
\end{example}

\begin{proposition}\label{p2.2}
The coring extension from example \ref{e2.1} satisfies the $(*)$-condition.
\end{proposition}
\begin{proof}
Note that for $(N,\rho_N)\in\Mm^\Dd$ and $x=n\otimes \sum\limits_{i\in I}\sigma_i(d_i)\in N\otimes \Cc$, by the definitions in Example \ref{e2.1} we have $\omega_{N,\Cc}(x)=n_0\otimes n_1\otimes \sum\limits_i\sigma_i(d_i)-\sum\limits_in\otimes d_{i1}\otimes \sigma_i(d_{i2})=\sum\limits_i(n_0\otimes n_1\otimes\sigma_i(d_i)-n\otimes d_{i1}\otimes\sigma_i(d_{i2}))$. If we denote by $D_i=N\otimes \Dd\otimes\sigma_i(\Dd)$ we have $N\otimes\Dd\otimes\Cc=\bigoplus\limits_iD_i=(N\otimes\Dd\otimes\Dd)^{(I)}$ and as $n_0\otimes n_1\otimes\sigma_i(d_i)-n\otimes d_{i1}\otimes\sigma_i(d_{i2})\in D_i$ we can see that $\omega_{N,\Cc}=\bigoplus\limits_{i\in I}\omega_{N,\Dd}$. Therefore it is enough to show that $\omega_{N,\Dd}$ is $\Cc$-pure. But $\Cc=\Dd^{(I)}$ as left $A$-modules, so it is enough to prove that $\omega_{N,\Dd}$ is $\Dd$-pure (because the tensor products commute with direct sums). It is easy to see that the sequence
$$0\longrightarrow N\stackrel{\rho_N}{\longrightarrow}N\otimes \Dd\stackrel{\omega_{N,\Dd}}{\longrightarrow}N\otimes\Dd\otimes\Dd$$
is exact, for $\omega_{N,\Dd}\circ\rho_N=0$ and if $\sum\limits_in_i\otimes d_i\in\Ker(\omega_{N,\Dd})$ then $\sum\limits_in_i\otimes d_i=\sum\limits_in_i\otimes d_{i1}\varepsilon_\Dd(d_i2)= \sum\limits_in_{i0}\otimes n_{i1}\varepsilon_\Dd(d_i)=\sum\limits_i(n_i\varepsilon_\Dd(d_i))_0\otimes (n_i\varepsilon_\Dd(d_i))_1\in\rho_N(N)$, thus $\rho_N(N)=\Ker\omega_{N,\Dd}$. If $x=\sum\limits_in_i\otimes d_i\otimes e_i\in \Ker(\omega_{N,\Dd}\otimes \Dd)$ then 
\begin{eqnarray*}
\sum\limits_in_i\otimes d_i\otimes e_i & = & \sum\limits_in_i\otimes d_{i1}\varepsilon_\Dd(d_{i2})\otimes e_i \\
& = & \sum\limits_in_{i0}\otimes n_{i1}\varepsilon_\Dd(d_i)\otimes e_i \;\;\;({\rm because\,}\omega_{N,\Dd}(x)=0)\\
& = & \sum\limits_i(n_i\varepsilon_\Dd(d_i))_0\otimes(n_i\varepsilon_\Dd(d_i))_1\otimes e_i
\end{eqnarray*}
showing that $x\in\im(\rho_N\otimes\Dd)$, so $\Ker(\omega_{N,\Dd}\otimes\Dd)=\im(\rho_N\otimes\Dd)$. Also if $n\otimes c\in\Ker(\rho_N\otimes\Dd)$ we have $n_0\otimes n_1\otimes d=0$ so $0=n_0\otimes\varepsilon_\Dd(n_1)d=n\otimes d$ and therefore $\rho_N\otimes\Dd$ is injective. Hence we get that the sequence
$$0\longrightarrow N\stackrel{\rho_N\otimes\Dd}{\longrightarrow}N\otimes \Dd\otimes\Dd\stackrel{\omega_{N,\Dd}\otimes\Dd}{\longrightarrow}N\otimes\Dd\otimes\Dd\otimes\Dd$$
is exact, and the proof is finished (in fact, it all follows as $N$ is a direct summand in $N\otimes \Dd$ as right $A$-modules).

\end{proof}

\begin{example}\label{e2.2}
Let $\Cc$ be an $A$-coring and $\Dd=A$ with the canonical Sweedler $A$-coring structure given by the comultiplication $A\ni a\mapsto1\otimes a=a\otimes 1\in A\otimes_AA$ and counit $A\ni a\mapsto a\in A$. Put $\lambda=\varepsilon_\Cc$. Then $\Cc\stackrel{\lambda}{\longrightarrow}A=\Cc\stackrel{\varepsilon_\Cc}{\longrightarrow}A$ is an extension of corings because $\lambda(c_1)\otimes\lambda(c_2)=\varepsilon_\Cc(c_1)\otimes_A\varepsilon_\Cc(c_2)=\varepsilon_\Cc(c_1)\varepsilon_\Cc(c_2)\otimes_A 1=\varepsilon_\Cc(c_1\varepsilon_\Cc(c_2))\otimes 1=\varepsilon_\Cc(c)\otimes 1=\lambda(c)_1\otimes\lambda(c)_{2}$. Moreover, the forgetful functor $\Mm^\Cc\longrightarrow \Mm_A$ associating to each right $\Cc$-comodule $M$ the underlying $A$-module $M$ is coincides with the corestriction functor associated to the morphism $\lambda=\varepsilon_\Cc$.
\end{example}

The following Proposition follows from \cite{BW}, Theorem 27.8. However we can obtain this from our generalization on Frobenius extensions of corings. 

\begin{proposition}{\bf (Frobenius corings)}
Let $\Cc$ be a coring. The extension of corings $\Cc\stackrel{\varepsilon_\Cc}{\longrightarrow}A$ from Example \ref{e2.2} is (right or left) Frobenius if and only if there is a Frobenius system $(e,\pi)$, $e\in\Cc^A=\{x\in\Cc\mid ax=xa,\,\forall\,a\in A\}$ and $\pi\in{}^\Cc\Hom^\Cc(\Cc\otimes_A\Cc,\Cc)$ such that
$$\pi(c\otimes e)=\pi(e\otimes c)=c$$
\end{proposition}
\begin{proof}
It is not difficult to see that we have an exact sequence 
$$0\longrightarrow \Cc\otimes\Cc\stackrel{1_{\Cc\otimes\Cc}}{\longrightarrow}\Cc\otimes\Cc\stackrel{\omega_{\Cc,\Cc}}{\longrightarrow}\Cc\otimes\Cc\otimes\Cc$$
because $\omega_{\Cc,\Cc}(c\otimes c')=c\otimes 1\otimes c'-c\otimes 1\otimes c'=0$ by the canonical $A$-comodule structure of $\Cc$. Therefore $\omega_{\Cc,\Cc}$ is always pure. Also for a right $A$-comodule (i.e. a right $A$-module) we can easily again see that $\omega_{N,A}=0$, and therefore the $(*)$-condition holds. Now if $\Cc\stackrel{\varepsilon_\Cc}{\longrightarrow}A$ is Frobenius, by Theorem \ref{th.Frob} we get $\alpha:A\longrightarrow \Cc$ an $A$-bimodule morphism and $\beta\in{}^\Cc\Hom^\Cc(\Cc\otimes\Cc,\Cc)$ such that equation (\ref{15}) holds. Put $\pi=\beta$ and $e=\alpha(1)$; then $ae=a\alpha(1)=\alpha(a)=\alpha(1)a=ea$ and $\alpha(a)=a\alpha(1)=ae=ea$. We get $\pi(e\otimes c)=\pi(e\otimes \varepsilon_\Cc(c_1)c_2)=\pi(e\varepsilon_\Cc(c_1)\otimes_A c_2)=\beta(\alpha(\varepsilon_\Cc(c_1))\otimes c_2)=c$ by (\ref{15}) and similarly $\pi(c\otimes e)=e$. Conversely if these conditions hold, define $\alpha=(a\mapsto ae)$ and $\beta=\pi$; then $\alpha\in{}^A\Hom^A(A,\Cc)$ and by the same computation as above we get that equation (\ref{15}) holds. Therefore by Theorem \ref{th.Frob} the extension $\Cc\stackrel{\varepsilon_\Cc}{\longrightarrow}A$ is a Frobenius extension.
\end{proof}

Following \cite{BW}, if an $A$-coring $\Cc$ satisfies the condition in the above proposition, $\Cc$ is said to be a {\it Frobenius coring}. By the above Proposition we see that this is a left-right symmetric concept. 

Let $\varphi:A\longrightarrow B$ be a morphism of rings. Recall that $A\stackrel{\varphi}{\longrightarrow}B$ is a Frobenius extension if $B$ is a finitely generated projective left (equivalently right) $A$-module and $B\simeq \Hom_A(B,A)$. This is equivalent to the fact that the induced forgetful functor ${}_B\Mm\longrightarrow {}_A\Mm$ is Frobenius, or equivalently, there are $E:B\longrightarrow A$ an $A$-bimodule morphism and an element $h=\sum\limits_{i}h_i\otimes g_i\in B\otimes_AB$ such that for all $b\in B$ we have
$$\sum_iE(bh_i)g_i\,=\,\sum_ih_iE(g_ib)\,=\,b$$
(we refer the reader to \cite{K} for these equivalent conditions). In this case the element $h$ is $B$-invariant, that is, $bh=hb$, $\forall\,b\in B$. Then it is easy to see that the existence of $h$ is equivalent to the existence of $u\in_B\Hom_B(B,B\otimes B)$, where $u(b)=hb$ (and conversely, $h=u(1)$). This comes from the isomorphism ${}_B\Hom_B(B,B\otimes_AB)\simeq \{h\in B\otimes_AB\mid bh=hb,\,\forall\,b\in B\}$, $u\mapsto u(1)$ (see \cite{BW}, section 27, and \cite{CMZ}). Therefore we can equivalently restate this as

\begin{proposition}\label{2.gen}
The ring extension $A\stackrel{\varphi}{\longrightarrow}B$ is Frobenius if and only if there are $E\in{}_A\Hom_A(B,A)$ and $u\in{}_B\Hom_B(B,B\otimes B)$ such that 
\begin{eqnarray}\label{20.a}
\mu\circ(\varphi E\otimes 1_B)\circ u\,=\,\mu\circ(1_B\otimes\varphi E)\circ u\,=\,1_B
\end{eqnarray}
where $\mu:B\otimes_A B\longrightarrow B$ is the multiplication of $B$ induced to $B\otimes_A B$.
\end{proposition}

\begin{remark}
By (\cite{BW}, 25.1, 27.7), if the extension $\varphi:A\longrightarrow B$ is Frobenius then $B\otimes_A B$ with the canonical Sweedler $B$-coring structure given by comultiplication $(b\otimes_A b')\mapsto (b\otimes_A 1)\otimes_B(1\otimes_A b')$ and counit $b\otimes_A b'\mapsto bb'$ is a Frobenius coring. The converse also holds  provided that $B$ is faithfully flat as left or right $A$-module. 
\end{remark}

In the case some restrictions are imposed on the base ring $A$, the restrictive conditions $(*)$ and "$\omega_{\Cc,\Cc}$ is pure in ${}_A\Mm$" can be eliminated. In particular, for extensions of coalgebras a theorem dual to the theorem characterizing Frobenius extensions of algebras (or rings, Proposition \ref{2.gen}) can be obtained.

\begin{example}{\bf Extensions of coalgebras}
An extension of corings $C\stackrel{\lambda}{\longrightarrow}D$ with $C,D$ coalgebras over a commutative ring $A$ will be called an extension of coalgebras.
\end{example}

\begin{theorem}\label{th.gen}
Let $A$ be a von Neumann regular ring (VNR). Then an extension of $A$-corings $\Cc\stackrel{\lambda}{\longrightarrow}\Dd$ is (left or right) Frobenius if and only if there are $\alpha\in{}^\Dd\Hom^\Dd(\Dd,\Cc)$ and $\beta\in{}^\Cc\Hom^\Cc(\Cc\square_{\Dd}\Cc,\Cc)$ such that 
$$
\beta(c_1\otimes\alpha\lambda(c_2)) \, = \, \beta(\alpha\lambda(c_1)\otimes c_2) \, = c
$$
for all $c\in \Cc$, equivalently,
\begin{eqnarray}\label{20.b}
\beta\circ(\alpha\lambda\otimes 1_\Cc)\circ\Delta\,=\,\beta\circ(1_\Cc\otimes \alpha\lambda)\circ\Delta\,=\,1_\Cc
\end{eqnarray}
In particular this holds for extensions of coalgebras over fields.
\end{theorem} 
\begin{proof}
If the base ring is VNR then all left and right modules are flat, and therefore the above mentioned conditions can be deleted from Theorem \ref{th.Frob}. Therefore, by the symmetry of the equation (\ref{15}) the theorem becomes left-right symmetric too.
\end{proof}

\begin{remark}
Theorem \ref{th.gen} gives a characterization of Frobenius extension of corings over VNR-rings (in particular of coalgebras over fields) which is completely dual to the characterization of Frobenius extensions of rings (and in particular of algebras over fields).
\end{remark}

\section{Finiteness Theorems}

The following proposition investigates when the coproduct of comodules on the category of $\Dd$-comodules indexed by a set $I$, $U=\bigoplus\limits_{i\in I}$, is a Frobenius functor. As the coproduct functor is a left adjoint to the diagonal functor $\delta:\Mm^\Dd\longrightarrow (\Mm^\Dd)^I$, this is equivalent to the fact that $U=\bigoplus$ is also a right adjoint to $\delta$, that is, it is also the product (of families indexed by $I$) in the category $\Mm^\Dd$, and the product and coproduct are isomorphic.

\begin{proposition}
The the extension of $A$-corings $\Cc\stackrel{\lambda}{\longrightarrow}\Dd$ from Example \ref{e2.1} with $\Dd\neq 0$ is a Frobenius extension (left or right) if and only if the set $I$ is finite. Consequently the coproduct functor $\bigoplus\limits_I$ on $\Mm^\Cc$ is Frobenius (and coincides with the product of comodules) if and only if $I$ is finite.
\end{proposition}
\begin{proof}
The statement follows from a result from \cite{I1}, namely Theorem 1.4, which shows that for a preadditive (and even for a more general type of) category if the coproduct (or product, or equivalently the diagonal functor) indexed by a set $I$ is a Frobenius functor then $I$ is finite (provided such a coproduct exists). However we can also see this from the results in the present paper. If $I$ is finite then it is easy to see that $U$ is also the product of comodules, with the projections being the canonical projections of the product of modules. To prove the converse, we first note that the $(*)$ condition holds by Proposition \ref{p2.2}. Then we can find $\alpha$ and $\beta$ as in Theorem \ref{th.Frob}. With the notations of Example \ref{e2.1}, let $\alpha_i=p_i\circ\alpha$; then $\alpha_i$ is a morphism in ${}^\Dd\Mm^\Dd$ because $p_i$ is too. First note that because $\beta$ is a morphism in $\Mm^\Cc$, for $d\in\Dd$ and $i,j\in I$ we have 
\begin{eqnarray*}
\beta(\sigma_j(d_1)\otimes\sigma_i(d_2)) & = & \varepsilon_\Cc(\beta(\sigma_j(d_1)\otimes\sigma_i(d_2))_1)\beta(\sigma_j(d_1)\otimes\sigma_i(d_2))_2\\
 & = & \varepsilon_\Cc(\beta(\sigma_j(d_1)\otimes \sigma_i(d_2)))\sigma_i(d_3)\in\sigma_i(\Dd)\\
\end{eqnarray*}
and similarly, as $\beta$ is a morphism in ${}^\Cc\Mm$ we get $\beta(\sigma_j(d_1)\otimes\sigma_i(d_2))\in\sigma_j(D)$ and therefore for $i\neq j$ we obviously get $\beta(\sigma_j(d_1)\otimes\sigma_i(d_2))=0$. Then for $d\in\Dd$ as $\alpha(d)=\sum\limits_j\sigma_j\alpha_j(d)$ we get 
\begin{eqnarray*}
d & = & \beta(\alpha\lambda(\sigma_i(d)_1)\otimes\sigma_i(d)_2) \;\;\; ({\rm by\,}(\ref{15}))\\
  & = & \beta(\alpha\lambda(\sigma_i(d_1))\otimes\sigma_i(d_2)) \;\;\; ({\rm definition\,of\,}\Delta_\Cc)\\
  & = & \beta(\alpha(d_1)\otimes \sigma_i(d_2)) \;\;\; ({\rm definition\,of\,}\lambda)\\
  & = & \sum\limits_{j\in I}\beta(\sigma_j\alpha_j(d_1)\otimes\sigma_i(d_2))\\
  & = & \sum\limits_{j\in I}\beta(\sigma_j(\alpha_j(d)_1)\otimes\sigma_i(\alpha_j(d)_2))\;\;\; ({\rm as\,}\alpha_j{\rm\,is\,a\,morphism\,in\,}\Mm^\Dd)\\
  & = & \sum\limits_{j\in I}\delta_{ij}\beta(\sigma_j(\alpha_j(d)_1)\otimes\sigma_i(\alpha_j(d)_2))\;\;\;({\rm where\,}\delta_{ij}\,{\rm is\,the\,Kroneker\,symbol})\\
  & = & \beta(\sigma_i(\alpha_i(d)_1)\otimes\sigma_i(\alpha_i(d)_2))
\end{eqnarray*}
This last equality obviously shows that for all $i\in I$, $\alpha_i$ is injective. But then for $d\in\Dd, d\neq 0$, we have $\alpha_i(d)\neq 0,\,\forall\,i\in I$ and therefore $I$ must be finite because the family $(\alpha_i(d))_{i\in I}\in\Dd^{(I)}$ is of finite support.
\end{proof}

For the coring extension of Example \ref{e2.2}, by \cite{BW}, 27.9 we have

\begin{proposition}
If $\Cc$ is a Frobenius coring (that is, the extension of corings from example \ref{e2.2} is Frobenius) then $\Cc$ is finitely generated projective as left and right $A$-module.
\end{proposition}

For extensions of coalgebras we can prove several interesting results parallel to existing ones for the extensions of algebras, namely several finiteness properties. We first prove a general finiteness theorem for Frobenius extensions of coalgebras. In what follows the ring $A$ is a field $K$ and the tensor product is always considered over $K$ unless otherwise specified.

\begin{theorem}\label{th.fin}{\bf The Finiteness of a Frobenius extension of coalgebras}
Let $C,D$ be two coalgebras over a commutative field $K$ and $C\stackrel{\lambda}{\longrightarrow}D$ a Frobenius extension of coalgebras. Then ${\rm dim}(C)\leq {\rm dim}(D)$ or they are both finite dimensional. 
\end{theorem}
\begin{proof} 
Take $\alpha$ and $\beta$ as in Theorem \ref{th.gen}. Then for $c\in C$ we have $\beta(\alpha\lambda(c_1)\otimes c_2)=c$. Let $(d_k)_{k\in\Lambda}$ be a $K$-basis for $D$ and for each $\alpha(d_k)$ choose an expression of $\Delta_C(\alpha(d_k))$ of the form $\Delta_C(\alpha(d_k))=\alpha(d_k)_1\otimes\alpha(d_k)_2=\sum\limits_lu_{kl}\otimes v_{kl}\in C\otimes C$. We show that $C$ is generated by the $\alpha(d_k)_1$'s, i.e. by the family of elements $(u_{kl})_{k,l}$; as for each $k\in\Lambda$ there are only a finite number of elements $u_{kl}$, the conclusion will follow. Denote by $C'=<(u_{kl})_{k,l}>$ the subspace of $C$ spanned by the family $(u_{kl})_{k,l}$. Take $c\in C$ and write $\Delta_C(c)=\sum\limits_ic_i\otimes e_i$. As $\lambda(c_i)\in D$, we have $\lambda(c_i)=\sum\limits_kd_kr_{ik}$, so $\alpha\lambda(c_i)=\sum\limits_k\alpha(d_k)r_{ik}$, $r_{ik}\in K$. Therefore we have $\alpha\lambda(c_1)_1\otimes \alpha\lambda(c_1)_2\otimes c_2=\sum\limits_i\alpha\lambda(c_i)_1\otimes \alpha\lambda(c_i)_2\otimes e_i=\sum\limits_i\sum\limits_k\alpha(d_k)_1\otimes\alpha(d_k)_2r_{ik}\otimes e_i=\sum\limits_{i,k,l}u_{kl}\otimes v_{kl}r_{ik}\otimes e_i$ so 
\begin{equation}\label{21}
\alpha\lambda(c_1)_1\otimes \alpha\lambda(c_1)_2\otimes c_2\,=\,\sum\limits_{i,k,l}u_{kl}\otimes v_{kl}r_{ik}\otimes e_i
\end{equation}
Note that we have an isomorphism $C\otimes(C\square_DC)\simeq (C\otimes C)\square_DC$ because we have an exact sequence
$$0\longrightarrow C\otimes(C\square_DC)\longrightarrow C\otimes C\otimes C\stackrel{C\otimes\omega_{C,C}}{\longrightarrow} C\otimes C\otimes C\otimes C$$
Then we can write $\alpha\lambda(c_1)_1\otimes \alpha\lambda(c_1)_2\otimes c_2=\sum\limits_sx_s\otimes T_s\in C\otimes(C\square_DC)$ with $T_s\in C\square_DC$. By a standard linear algebra argument, we can take $(T_s)_s$ to be linearly independent (just take an expression of $\sum\limits_sx_s\otimes T_s$ with a minimal number of tensor monomials of the type $x\otimes T$, $x\in C$, $T\in C\square_DC$). Note that $\alpha\lambda(c_1)_1\otimes \beta(\alpha\lambda(c_1)_2\otimes c_2)=(C\otimes \beta)\circ(\Delta_C\otimes C)\circ(\alpha\lambda\otimes C)\circ\Delta_\Cc$. Therefore we have 
\begin{eqnarray*}
\sum\limits_sx_s\varepsilon_C\beta(T_s) & = & \alpha\lambda(c_1)_1\varepsilon_C\beta(\alpha\lambda(c_1)_2\otimes c_2) \\
 & = & \beta(\alpha\lambda(c_1)\otimes c_2)_1\varepsilon_C(\beta(\alpha\lambda(c_1)\otimes c_2)_2) \;\;\;({\rm because\,}\beta\,{\rm is\,a\,morphism\,in\,}{}^C\Mm)\\
 & = & \beta(\alpha\lambda(c_1)\otimes c_2) \;\;\;({\rm by\,the\,counit\,property})\\
 & = & c \;\;\;({\rm by\,}(\ref{15}))
\end{eqnarray*}
so we get 
\begin{equation}
c\,=\,\sum\limits_sx_s\varepsilon_C\beta(T_s) \label{22}
\end{equation}
As $(T_s)_s$ are independent there are $T_s^*\in(C\square_DC)^*$ such that $T_p^*(T_s)=\delta_{ps}$; as $C\square_DC\subseteq C\otimes C$ we can find $U_s^*\in(C\otimes C)^*$ such that $U_{s}^*\vert_{ C\square_DC}=T_s^*$. Then by (\ref{21}) we have $\sum\limits_sx_s\otimes T_s=\sum\limits_{i,k,l}u_{kl}\otimes v_{kl}r_{ik}\otimes e_i$, so we get $x_p=\sum\limits_sx_sT_p^*(T_s)=\sum\limits_sx_sU_p^*(T_s)=\sum\limits_{i,k,l}u_{kl}U_p^*(r_{ik}v_{kl}\otimes e_i)\in C'$ and therefore by (\ref{22}) we find that $c=\sum\limits_s\varepsilon_C\beta(T_s)x_s\in C'$. Thus $C'=C$ and the proof is finished.
\end{proof}

\begin{proposition}
Let $C\stackrel{\lambda}{\longrightarrow}D$ be an extension of coalgebras (over a field $K$) and denote $E=\lambda(C)$. Then  $C\square_DC=C\square_EC$ and $C^*\otimes_{D^*}C^*=C^*\otimes_{E^*}C^*$.
\end{proposition}
\begin{proof}
We have that $E$ is a subcoalgebra of $D$. Write $\lambda=j\circ i$ with $i:C\longrightarrow E$ being the corestriction of $\lambda$ and $j:E\subseteq D$ the canonical inclusion. Note that we have a commutative diagram
\begin{diagram}
0 & \rTo & C\square_DC & \rTo & C\otimes C & \rTo^{\omega_{C,C}^D} & C\otimes D\otimes C \\
  &			 &             &      &  \dEqual   &                       &  \uTo_{C\otimes i\otimes C} \\
0 & \rTo & C\square_EC & \rTo & C\otimes C & \rTo_{\omega_{C,C}^E} & C\otimes E\otimes C\\
\end{diagram}
which shows that $C\square_DC=C\square_EC$ because $C\otimes i\otimes C$ is injective (because $i$ is). Similarly, as there is a morphism $j^*:D^*\longrightarrow E^*$, we have an epimorphism $\varphi_1:C^*\otimes_{D^*}C^*\longrightarrow C^*\otimes_{E^*}C^*$ taking $g^*\otimes_{D^*}h^*$ to $g^*\otimes_{E^*}h^*$, because the application $C^*\otimes C^*\ni(g^*,h^*)\mapsto g^*\otimes_{E^*}h^*\in C^*\otimes_{E^*}C^*$ is $D^*$-balanced. But as for $e^*\in E^*$ and $g^*,h^*\in C^*$, there is $d^*\in D^*$ with $j^*(d^*)=e^*$, we have 
\begin{eqnarray*}
h^*\cdot e^*\otimes_{D^*} g^* & = & h^**i^*(e^*)\otimes_{D^*}g^* \,=\, h^**i^*j^*(d^*)\otimes_{D^*}g^* \\
 & = & h^**\lambda^*(d^*)\otimes_{D^*}g^* \,=\, h^*\cdot d^*\otimes_{D^*}g^* \\
 & = & h^*\otimes_{D^*}d^*\cdot g^* \,=\, h^*\otimes_{D^*}\lambda(d^*)*g^* \\
 & = & h^*\otimes_{D^*} i^*(e^*)*g^* \,=\, h^*\otimes_{D^*}e^*\cdot g^*
\end{eqnarray*}
This shows that the application $\varphi_2:C^*\otimes_{E^*}C^*\longrightarrow C^*\otimes_{D^*}C^*$, $\varphi_2(g^*\otimes_{E^*}h^*)=g^*\otimes_{D^*}h^*$ is well defined and it is obviously an inverse for $\varphi_1$. 
\end{proof}

\begin{lemma}\label{lemma}\label{izo.1}
Let $C\stackrel{\lambda}{\longrightarrow}D$ be an extension of coalgebras (over a field $K$) such that $C$ is finite dimensional. Then there is an isomorphism of $C^*$-bimodules $C^*\otimes_{D^*}C^*\simeq (C\square_DC)^*$ given by $g^*\otimes_{D^*} h^*\longmapsto(g\otimes h\mapsto g^*(g)h^*(h))$. 
\end{lemma}
\begin{proof}
We have an isomorphism of vector spaces $\varphi:C^*\otimes C^*\stackrel{\sim}{\longrightarrow} (C\otimes C)^*$ given by $\varphi(g^*\otimes h^*)=(g\otimes h\mapsto g^*(g)h^*(h))$. First note that by the previous proposition, replacing $D$ with $\lambda(C)$ we may assume that $D$ is finite dimensional too (the $C^*$-module structure is also preserved). Let $\eta:C\square_DC\longrightarrow C\otimes C$ be the inclusion morphism and $p:C^*\otimes C^*\longrightarrow C^*\otimes_{D^*}C^*$ the canonical morphism $p(g^*\otimes h^*)=g^*\otimes_{D^*}h^*$; define $\pi:C^*\otimes_{D^*}C^*\longrightarrow (C\square_DC)^*$ by $\pi(g^*\otimes_{D^*}h^*)=\eta^*\varphi(g^*\otimes h^*)$. Note that $\pi$ is well defined as for $d^*\in D^*$ and $g\otimes h\in C\square_DC$ we have $g_1\otimes\lambda(g_2)\otimes h=g\otimes\lambda(h_1)\otimes h_2$ so $\eta^*\varphi(g^*\cdot d^*\otimes h^*)(g\otimes h)=\eta^*\varphi(g^**\lambda^*(d^*)\otimes h^*)(g\otimes h)=g^*(g_1)d^*(\lambda(g_2))h^*(h)=g^*(g)d^*(\lambda(h_1))h^*(h_2)=\eta^*\varphi(g^*\otimes \lambda^*(d^*)*h^*)=\eta^*\varphi(g^*\otimes d^*\cdot h^*)$ ($*$ is the convolution product of $C^*$). Also note that $\pi$ is surjective as $\pi\circ p=\eta^*\circ\varphi$ and $\eta^*$ is surjective. As we have epimorphisms
\begin{diagram}
(C\otimes C)^*=C^*\otimes C^* & \rOnto^{p} & C^*\otimes_{D^*}C^* & \rOnto^{\pi} & (C\square_DC)^* 
\end{diagram}
we get monomorphisms and a commutative diagram
\begin{diagram}
(C\square_DC)^{**}          & \rInto^{\pi^*} & (C^*\otimes_{D^*}C^*)^* & \rInto^{p^*} & (C^*\otimes C^*)^*\simeq (C\otimes C)^{**} \\
\dEqual_{\Psi_{C\square_DC}}&                &                         &                 &  \dEqual^{\Psi_{C\otimes C}} \\
C\square_DC                 &                &   \rTo^{\varphi}        &                 &  C\otimes C
\end{diagram}
Here, for a vector space $V$ we use the identification between $V$ and $V^{**}$ given by the isomorphism $\Psi_V:V\longrightarrow V^{**}$, $\Psi_V(v)=(v^*\mapsto v^*(v))\in V^{**}$. We prove that $\pi^*$ is surjective and this will show that $\pi$ is also injective. Pick $\Lambda\in(C^*\otimes_{D^*}C^*)^*$. Then $p^*(\Lambda)=\Lambda\circ p\in (C^*\otimes C^*)^*\simeq C\otimes C$ so by the identification we made we can find $c\otimes e\in C\otimes C$ such that $p^*(\Lambda)(g^*\otimes h^*)=g^*(c)h^*(e)$, $\forall\,g^*\otimes h^*\in C^*\otimes C^*$. For $d^*\in D^*$ we have 
\begin{eqnarray*}
g^*(c_1)d^*(\lambda(c_2))h^*(e) & = & p^*(\Lambda)(g^**\lambda^*(d^*)\otimes h^*)(c\otimes e)\\
 & = & \Lambda(g^*\cdot d^*\otimes_{D^*} h^*)(c\otimes e) \\
 & = & \Lambda(g^*\otimes_{D^*} d^*\cdot h^*)(c\otimes e) \\
 & = & p^*(\Lambda)(g^*\otimes \lambda(d^*)*g^*)(c\otimes e) \\
 & = & g^*(c)d^*(\lambda(e_1))h^*(e_2)
\end{eqnarray*}
showing that $g^*(c_1)d^*(\lambda(c_2))h^*(e)\,=\,g^*(c)d^*(\lambda(e_1))h^*(e_2)$ for all $g^*,h^*\in C^*$ and $d^*\in D^*$ so $U(c_1\otimes \lambda(c_2)\otimes e)=U(c\otimes\lambda(e_1)\otimes e_2)$, $\forall\,U\in(C\otimes D\otimes C)^*$ and therefore $c_1\otimes \lambda(c_2)\otimes e=c\otimes \lambda(e_1)\otimes e_2$. Thus $c\otimes e\in C\square_DC$, that is $c\otimes e=\varphi(c\otimes e)$. Hence $p^*(\Lambda)=\Psi_{C\otimes C}(\varphi(c\otimes e))=p^*(\pi^*(\Psi_{C\square_DC}(c\otimes e)))$, so $\Lambda=\pi^*(\Psi_{C\square_DC}(c\otimes e))$. Thus $\pi^*$ is a bijection. \\
Finally, for an element $\chi=\sum\limits_ih_i^*\otimes_{D^*} g_i^*\in C^*\otimes_{D^*}C^*$ and $c^*\in C^*$, $c\otimes e\in C\square_DC$ we have $(c^*\pi(\chi))(c\otimes e)=\pi(\chi)((c\otimes e)\cdot c^*))=c^*(c_1)\pi(\chi)(c_2\otimes e)=\sum\limits_ic^*(c_1)h_i^*(c_2)g_i^*(e)=\pi(c^*\cdot h_i^*\otimes_{D^*}g_i^*)(c\otimes e)=\pi(c^*\cdot\chi)(c\otimes e)$. This shows that $\pi$ is a morphism of left (and similarly of right) $C^*$-modules.
\end{proof}

\begin{proposition}\label{3.inj}
If $C\stackrel{\lambda}{\longrightarrow}D$ is a Frobenius extension of coalgebras, then $C$ is injective as left and also as right $D$-comodule. 
\end{proposition}
\begin{proof}
By the definition of Frobenius extensions, we have that the functor $F=-\square_DC:\Mm^D\longrightarrow \Mm^C$ is Frobenius so it has the same left and right adjoint. As $\Mm^C$ and $\Mm^D$ are abelian categories, we get that $F$ is left and right exact. Therefore by \cite{DNR}, 2.4.23 we have that $C$ is injective as left $D$-comodule. By Theorem \ref{th.gen} we have the left-right symmetry of Frobenius extensions of coalgebras and therefore we also get that $C\square_D-:{}^D\Mm\longrightarrow {}^C\Mm$ is a Frobenius functor and that $C$ is also injective as right $D$-comodule.
\end{proof}

\begin{proposition}\label{izo.2}
Let $C\longrightarrow D$ be an extension of coalgebras with $C$ finite dimensional. Then the application $\Phi:{}^D\Hom^D(D,C)\ni\alpha\longrightarrow \Phi(\alpha)=\alpha^*\in{}_{D^*}\Hom_{D^*}(C^*,D^*)$ is well defined and bijective, where ${}^D\Hom^D$ (${}_{D^*}\Hom_{D^*}$) represents the set of morphisms of $D$-bicomodules (respectively $D^*$-bimodules).
\end{proposition}
\begin{proof}
It is not difficult to see that $\Phi$ is injective, as if $\alpha^*=0$, then $c^*(\alpha(d))=0$ for all $d\in D$ and $c^*\in C^*$ and therefore $\alpha(d)=0,\,\forall\,d\in D$. Let $u:C^*\longrightarrow D^*$ be a morphism of $D^*$-bimodules.  Then $I=\im(u)$ is an ideal of $D^*$ of finite dimension and thus $E=I^\perp=\{d\in D\mid d^*(d)=0,\,\forall\,d^*\in I\}$ is a subcoalgebra of $D$ which has finite codimension as $I$ is finite dimensional. We also have that $I$ is closed in the finite topology on $D^*$ (for example, by \cite{DNR}, Corollary 1.2.12), so $E^\perp=\{d^*\in D^*\mid d^*(d)=0,\,\forall\,d\in E\}=(I^\perp)^\perp=I$. Let $p:D\longrightarrow D/E$ be the canonical projection, $i:I\hookrightarrow D^*$ the inclusion morphism and $v:C^*\longrightarrow I$ the corestriction of $u$. There is an isomorphism $\theta:(D/E)^*\longrightarrow I=E^\perp$ taking $h\in(D/E)^*$ to $\theta(h)=h\circ p\in I$, and then we have a commutative diagram
\begin{diagram}
C^* & \rTo^{v} 							&       I                 & \stackrel{i}{\hookrightarrow} & D^* \\ 
    & \rdTo_{\overline{v}}  & \uTo_\simeq^{\theta}    & \ruTo_{p^*}                   &      \\
    &                       & (\frac{D}{E})^*         &                               &
\end{diagram}
where $\overline{v}=\theta^{-1}\circ v$. But $C$ and $D/E$ are finite dimensional vector spaces and therefore there is $q:D/E\longrightarrow C$ such that $\overline{v}=q^*$. Hence we obtain $u=i\circ v=i\circ \theta\circ \overline{v}=p^*\circ \overline{v}=p^*\circ q^*=(q\circ p)^*$, so $u=\alpha^*$ for $\alpha=q\circ p$. \\
Now the following sequence of equivalences shows that $\alpha$ is a $D$-bicomodule morphism if and only if $\alpha^*$ is a $D^*$-bimodule morphism, and thus $\Phi$ is well defined (denote by $*$ the convolution product of $C^*$):
\begin{eqnarray*}
  & & \alpha {\rm\,is\,a\,morphism\,of\,right\,}D-{\rm\,comodules} \;\;\; \Leftrightarrow\\
\alpha(d_1)\otimes d_2 & = & \alpha(d)_1\otimes\lambda(\alpha(d)_2),\,\forall\,d\in D \;\;\; \Leftrightarrow \\
c^*(\alpha(d_1))d^*(d_2)  & = & c^*(\alpha(d)_1)d^*(\lambda(\alpha(d)_2)),\,\forall\,d\in D,\,\forall\,c^*\in C^*,\,\forall\,d^*\in D^* \;\;\; \Leftrightarrow \\
(\alpha^*(c^*)*d^*)(d) & = & (c^**\lambda^*(d^*))(\alpha(d)),\,\forall\,d\in D,\,\forall\,c^*\in C^*,\,\forall\,d^*\in D^* \;\;\; \Leftrightarrow  \\
\alpha^*(c^*)*d^* & = & \alpha^*(c^*\cdot d^*),\,\forall\,c^*\in C^*,\,\forall\,d^*\in D^* \;\;\; \Leftrightarrow  \\
  & & \alpha^*\,{\rm\,is\,a\,morphism\,of\,right\,}D^*-{\rm\,modules} 
\end{eqnarray*}
(and similarly $\alpha$ is a morphism of left $D$-comodules if and only if $\alpha^*$ is a morphism of left $D^*$-modules.)
\end{proof}

\begin{proposition}\label{3.p1}
Let $C$ be a finite dimensional coalgebra and $C\stackrel{\lambda}{\longrightarrow}D$ an extension of coalgebras. Then the following assertions are equivalent:
\begin{itemize}
\item[(i)] $C\stackrel{\lambda}{\longrightarrow}D$ is a Frobenius extension of coalgebras.
\item[(ii)] $D^*\stackrel{\lambda^*}{\longrightarrow}C^*$ is a Frobenius extension of algebras (rings).
\end{itemize}
Moreover, in this case $C$ is finitely cogenerated and injective as left and also as right $D$-comodule and $C^*$ is finitely generated projective as right and also as left $D^*$-module.
\end{proposition}
\begin{proof}
(i)$\Rightarrow$(ii) By Theorem \ref{th.gen} we can find $\alpha,\beta$ such that (\ref{15}) and (\ref{20.b}) hold. Put $E=\alpha^*\in{}_{D^*}\Hom_{D^*}(C^*,D^*)$ by Proposition \ref{izo.2}. As $\beta:C\square_DC\longrightarrow C$ is a morphism of $C$-bicomodules, we can easily see that $\beta^*:C^*\longrightarrow (C\square_DC)^*$ is a morphism of $C^*$-bimodules, because the category of finite dimensional $C$-bicomodules is in duality with the category of $C^*$-bimodules. Then by Lemma \ref{izo.1} we may view $\beta^*$ as a morphism of $C^*$-bimodules $u=\beta^*:C^*\longrightarrow C^*\otimes_{D^*}C^*\simeq(C\square_DC)^*$. Now by dualizing equation (\ref{20.b}) we can easily see that we obtain (\ref{20.a}), where $\mu=\Delta^*$ and $1_{C^*}=\varepsilon_C$ are the multiplication and unit of $C^*$ and $\varphi=\lambda^*$, showing that $D^*\stackrel{\lambda^*}{\longrightarrow}C^*$ is a Frobenius extension.\\ 
(ii)$\Rightarrow$(i) By Proposition \ref{2.gen} we find $E\in{}_{D^*}\Hom_{D^*}(C^*,D^*)$ and $u\in{}_{C^*}\Hom_{C*}(C^*,C^*\otimes_{D^*}C^*)$ such that (\ref{20.a}) holds (with $\mu$ and $\varphi$ as above). Then by Lemma \ref{izo.1}, as $C^*\otimes_{D^*}C^*\simeq (C\square_DC)^*$ as $C^*$-bimodules and $(C\square_DC)^*$ and $C^*$ are finite dimensional, there is a morphism of $C$-bicomodules $\beta:C\square_DC\longrightarrow C$ such that $\beta^*=u$. Also, by Proposition \ref{izo.2}, $E=\alpha^*$ with $\alpha:D\longrightarrow C$ a morphism of $D$-bicomodules. Then by (\ref{20.a}), $c^*=(\Delta^*\circ(\lambda^*\alpha^*\otimes 1_{C^*})\circ \beta^*)(c^*)=(\Delta^*\circ(\alpha\lambda\otimes 1_C)^*)(c^*\circ\beta)=c^*\circ \beta\circ(\alpha\lambda\otimes 1_C)\circ \Delta$ for all $c^*\in C^*$ and therefore $\beta\circ(\alpha\lambda\otimes 1_C)\circ \Delta=1_C$ and similarly $\beta\circ(1_C\otimes\alpha\lambda)\circ \Delta=1_C$. Thus (\ref{20.b}) holds and $C\stackrel{\lambda}{\longrightarrow}D$ is a Frobenius extension of coalgebras. \\
Finally if (i) and (ii) hold, by Proposition \ref{3.inj} $C$ is injective in $\Mm^D$ (and also in ${}^D\Mm$), so there is a monomorphism of right (left) $C\hookrightarrow D^{(I)}$ for a set $I$. But as $C$ is finite dimensional, we may obviously assume that the set $I$ is finite. Thus we get $C\hookrightarrow D^n$ as right (or left) $D$-comodules for some $n\in\NN$, and in fact $C$ is a direct summand of $D^n$ as it is injective. Also as $D^*\stackrel{\lambda^*}{\longrightarrow}C^*$ is a Frobenius extension of rings we get that $C^*$ is finitely generated and projective as right (and also as left) $D^*$-module (or as $C$ splits off in some $D^n$). 
\end{proof}

\begin{proposition}
Assume $C\stackrel{\lambda}{\longrightarrow}D$ is an extension of coalgebras with $D$ finite dimensional. Then the following are equivalent:
\begin{itemize} 
\item[(i)] $C\stackrel{\lambda}{\longrightarrow}D$ is a Frobenius extension of coalgebras. 
\item[(ii)] $D^*\stackrel{\lambda^*}{\longrightarrow}C^*$ is a Frobenius extension of algebras. 
\end{itemize}
In this case $C$ is finite dimensional and injective finitely cogenerated as right (and also as left) $D$-comodule and $C^*$ is finitely generated as left (and also as right) $D^*$-module.
\end{proposition}
\begin{proof}
(i)$\Rightarrow$(ii) By Theorem \ref{th.fin} we get that $C$ is finite dimensional as $D$ is. Then $\Mm^C$ coincides with ${}_{C^*}\Mm$ and $\Mm^D$ with ${}_{D^*}\Mm$, as any $C^*$ ($D^*$)-module is in this case a $C$ ($D$)-comodule. Note that in this case the corestriction functor $U:\Mm^C\longrightarrow \Mm^D$ induced by $\lambda$ identifies with the forgetful functor $H:{}_{C^*}\Mm\longrightarrow {}_{D^*}\Mm$. Indeed if $(M,\rho_M)\in \Mm^C$, then $M\in \Mm^D$ by the coaction $m\mapsto m_0\otimes \lambda(m_1)$ where $m_0\otimes m_1=\rho_M(m)$. Then $M$ becomes a left $D^*$-module by $d^*\cdot m=m_0d^*(\lambda(m1))$ for $m\in M$ and $d^*\in D^*$. But then $d^*\cdot m=m_0\lambda^*(d^*)(m_1)=\lambda^*(d^*)\cdot m$ which shows that the $D^*$-module structure of $M$ coincides with the one obtained through the forgetful functor $H$. Therefore, if $C\stackrel{\lambda}{\longrightarrow}D$ is a Frobenius extension, $U$ is a Frobenius functor and therefore $H$ is Frobenius, proving (ii).\\
(ii)$\Rightarrow$(i) If $D^*\stackrel{\lambda^*}{\longrightarrow}C^*$ is Frobenius then $C^*$ must be finitely generated (and projective) as right (and also as left) $D^*$-module and as $D^*$ is finite dimensional, we get that $C$ is finite dimensional too. Then as in (i)$\Rightarrow$(ii) we obtain that the induced functor $U$ is Frobenius and therefore $C\stackrel{\lambda}{\longrightarrow}D$ is a Frobenius extension. \\
The last statement follows easily as in the proof of Proposition \ref{3.p1}.
\end{proof}

\begin{remark}
It can be seen that in a wide range of situations whenever a functor is Frobenius, a certain finiteness theorem holds. Then it is natural to ask whether a general theorem can be proved; that is, given a Frobenius pair $(F,G)$ between two abelian categories (or Grothendieck, or with other additional properties) is there a finiteness property that holds for this general context and which generalizes all these theorems regarding Frobenius extensions of (co)rings?
\end{remark}

\bigskip\bigskip\bigskip

\begin{center}
\sc Acknowledgment
\end{center}
The author wishes to thank his Ph.D. adviser C. N\u ast\u asescu for very useful remarks on the subject as well as for his continuous support throughout the past years. He would also like to address special thanks to the referee for a very helpful report.


\vspace*{3mm} 
\begin{flushright}
\begin{minipage}{148mm}\sc\footnotesize

Miodrag Cristian Iovanov\\
University of Bucharest, Faculty of Mathematics, Str. Academiei 14\\ 
RO-010014, Bucharest, Romania\\
{\it E--mail address}: {\tt
yovanov@walla.com}\vspace*{3mm}

\end{minipage}
\end{flushright}
\end{document}